\documentclass{amsart}
\usepackage{amssymb}
\usepackage{amscd}
\usepackage{eucal}
\usepackage{latexsym}

\newcommand{\jj}{\vee}
\newcommand{\mm}{\wedge}
\newcommand{\JJ}{\bigvee}
\newcommand{\MM}{\bigwedge}
\newcommand{\JJm}[2]{\JJ(\,#1\mid#2\,)}
\newcommand{\MMm}[2]{\MM(\,#1\mid#2\,)}

\newcommand{\uu}{\cup}

\newcommand{\ci}{\subseteq}
\newcommand{\nin}{\notin}
\newcommand{\es}{\varnothing}
\newcommand{\set}[1]{\{#1\}}
\newcommand{\setm}[2]{\{\,#1\mid#2\,\}}
\def\vv<#1>{\langle#1\rangle}

\newcommand{\ga}{\alpha}
\newcommand{\gb}{\beta}

\renewcommand{\ge}{\varepsilon}
\newcommand{\gf}{\varphi}
\renewcommand{\gg}{\gamma}
\newcommand{\gh}{\eta}
\newcommand{\gi}{\iota}

\newcommand{\gk}{\kappa}

\newcommand{\gm}{\mu}
\newcommand{\gn}{\nu}
\newcommand{\go}{\omega}
\newcommand{\gp}{\pi}

\newcommand{\gr}{\varrho}
\newcommand{\gs}{\sigma}

\newcommand{\gx}{\xi}
\newcommand{\gy}{\psi}
\newcommand{\gz}{\zeta}

\newcommand{\gF}{\Phi}

\newcommand{\gQ}{\Theta}

\newcommand{\tbf}{\textbf}
\newcommand{\tup}{\textup}

\newcommand{\E}[1]{\mathcal{#1}}

\providecommand{\bysame}{\makebox[3em]{\hrulefill}\thinspace}
\newcommand{\q}{\quad}
\newcommand{\iso}{\cong}
\def\con#1=#2(#3){#1\equiv#2\pod{#3}}

\theoremstyle{plain}
\newtheorem{lemma}{Lemma}[section]
\newtheorem{theorem}[lemma]{Theorem}
\newtheorem{proposition}[lemma]{Proposition}
\newtheorem{corollary}[lemma]{Corollary}
\newtheorem{claim}{Claim}

\newtheorem*{tha}{Theorem 1}
\newtheorem*{thb}{Theorem 2}
\newtheorem*{thc}{Theorem 3}

\theoremstyle{definition}
\newtheorem{definition}[lemma]{Definition}

\newtheorem{problem}{Problem}
\newtheorem*{remark}{Remark}

\numberwithin{equation}{section}

\newcommand{\qedc}{{\qed}~{\rm Claim~{\theclaim}.}}
\newcommand{\sqedc}{{\qed}~{\rm Claim.}}

\newenvironment{cproof} {\begin{proof}[Proof of Claim.]}
{\qedc\renewcommand{\qed}{}\end{proof}}

\DeclareMathOperator{\Con}{Con}
\DeclareMathOperator{\Conc}{Con_c}
\DeclareMathOperator{\Res}{Res}

\newcommand{\jz}{$\set{\jj,0}$}

\newcommand{\zero}{\mathbf{0}}
\newcommand{\one}{\mathbf{1}}
\newcommand{\two}{\mathbf{2}}

\newcommand{\tgr}{\tilde{\gr}}
\newcommand{\tga}{\tilde{\ga}}
\newcommand{\tgb}{\tilde{\gb}}
\newcommand{\bga}{\bar{\ga}}
\newcommand{\bgb}{\bar{\gb}}

\newcommand{\jhom}{join-ho\-mo\-mor\-phism}
\newcommand{\mhom}{meet-ho\-mo\-mor\-phism}
\newcommand{\cjh}{complete join-ho\-mo\-mor\-phism}
\newcommand{\cmh}{complete meet-ho\-mo\-mor\-phism}

\newcommand{\bgi}{\boldsymbol{\gi}}

\newcommand{\Dc}{{\E{D}_{\mathrm{c}}}}
\newcommand{\Dac}{{\E{D}_{\mathrm{ac}}}}

\renewcommand{\AA}{\mathbf{A}}
\newcommand{\BB}{\mathbf{B}}
\newcommand{\FF}{\mathbf{F}}
\newcommand{\VV}{\E{V}}
\newcommand{\LL}{$\E{L}$}
\newcommand{\Cj}{\mathbf{C}_{\jj}}
\newcommand{\Cm}{\mathbf{C}_{\mm}}
\newcommand{\FV}{\mathrm{F}_{\VV}}
\newcommand{\id}{\mathrm{id}}
\newcommand{\res}{\mathbin{\restriction}}

\newcommand{\URPI}{$\mathrm{URP}_1$}
\newcommand{\URPIM}{$\mathrm{URP}_1^-$}

\newcommand{\tvi}{\vrule height 12pt depth 6pt width 0pt}

\newcommand{\cc}[1]{\hfill\kern.7em#1\kern.7em\hfill}

\begin{document}

\title[Representations by lattices with permutable
congruences]%
{Simultaneous representations of semilattices\\
by lattices with permutable congruences}

 \author[J.~T\r uma]{Ji\v r\'\i\ T\r uma}
 \address{Department of Algebra\\
          Faculty of Mathematics and Physics\\
          Sokolovsk\'a 83\\
          Charles University\\
          186 00 Praha 8\\
          Czech Republic}
 \email{tuma@karlin.mff.cuni.cz}
 \thanks{The first author was partially supported by
 GA \v CR 201/97/1162.}

 \author[F.~Wehrung]{Friedrich Wehrung}
 \address{C.N.R.S., E.S.A. 6081\\
          Universit\'e de Caen, Campus II\\
          D\'epartement de Math\'ematiques\\
          B.P. 5186\\
          14032 CAEN Cedex\\
          FRANCE}
 \email{wehrung@math.unicaen.fr}
 \urladdr{http://www.math.unicaen.fr/\~{}wehrung}

\keywords{Lattice, semilattice, distributive,
congruence, permutable, Uniform Refinement Property,
congruence-splitting, diagram of lattices,
congruence representation}
\subjclass{06A12, 06B10}

\date{\today}

\begin{abstract}
The Congruence Lattice Problem (CLP), stated by R.~P. Dilworth
in the forties, asks whether every distributive \jz-semilattice
$S$ is isomorphic to the semilattice $\Conc L$ of
compact congruences of a lattice $L$.

While this problem is still
open, many partial solutions have been obtained, positive and
negative as well. The solution to CLP is known to be positive
for all $S$ such that $|S|\leq\aleph_1$. Furthermore,
one can then take $L$ \emph{with permutable congruences}. This
contrasts with the case where $|S|\geq\aleph_2$, where there are
counterexamples $S$ for which $L$ \emph{cannot} be, for example,
sectionally complemented.
We prove in this paper that the lattices of these counterexamples
cannot have permutable congruences as well.

We also isolate \emph{finite, combinatorial
analogues} of these results. All the ``finite'' statements
that we obtain are \emph{amalgamation properties} of the
$\Conc$ functor. The strongest known
positive results, which originate in earlier work by the first
author, imply that many diagrams of semilattices indexed by
the \emph{square} $\two^2$ can be lifted with respect to the
$\Conc$ functor.

We prove that the latter results cannot be extended to the
\emph{cube}, $\two^3$. In particular, we give an example of a
cube diagram of finite Boolean semilattices and semilattice
embeddings that cannot be lifted, with respect to the
$\Conc$ functor, by lattices with permutable congruences.

We also extend many of our results to lattices with \emph{almost
permutable congruences}, that is, $\ga\jj\gb=\ga\gb\uu\gb\ga$,
for all congruences $\ga$ and $\gb$.

We conclude the paper with a very short proof that no functor
from finite Boolean semilattices to lattices can lift the
$\Conc$ functor on finite Boolean semilattices.

\end{abstract}

\maketitle

\section*{Introduction}

The classical Congruence Lattice Problem, posed by R.~P.
Dilworth in the early forties, formulated in the
language of semilattices, asks whether any distributive
\jz-semilattice $S$ is isomorphic to the
semilattice of compact congruences $\Conc L$ of a lattice
$L$. This is probably the most well-known open problem in
lattice theory, see \cite{GrSc} for a survey.

The answer to CLP is known to be positive if $|S|\leq\aleph_1$.
In fact, in that case, one can take $L$ relatively
complemented, locally finite, with zero, see
\cite{GLWe}; this result can also be derived from the results
of \cite{Tuma}. Still in case $|S|\leq\aleph_1$, one can take
$L$ sectionally complemented, \emph{modular} (see \cite{Wehr3}).
It is an important observation that in both cases, \emph{any two
congruences of $L$ permute}, see \cite{Dilw}.

Partial negative answers to this problem are obtained by
applying a universal construction, described in
the second author's paper
\cite{Wehr1}, that yields, for every infinite cardinal
number $\gk$, a ``complicated'' distributive \jz-semilattice
$S_\gk$ of cardinality
$\gk$. A variant of $S_\gk$, with similar properties, is
constructed in \cite{PTu}. Although its precise construction is
relatively complicated, this variant could be described as a
distributive
\jz-semilattice ``freely generated'' by elements
$a_\gx$, $b_\gx$, for $\gx<\gk$, subject to the relation
$a_\gx\jj b_\gx=\text{constant}$,
for all $\gx$. If $\gk\geq\aleph_2$, then $S_\gk$ is not
isomorphic to $\Conc L$ if $L$ is, say, sectionally
complemented, or, more generally, \emph{congruence-splitting},
as defined in \cite{Wehr2}.

The key observation is that if $L$ is congruence-splitting,
then $\Conc L$ satisfies a certain infinite axiom of the
language of semilattices, called the \emph{Uniform
Refinement Property}, URP. The main results can be
summarized in the following theorem, see
\cite{PTWe,Wehr1,Wehr2}.

\goodbreak
\begin{tha}\hfill
\begin{enumerate}
\item The class of congruence-splitting lattices contains
the class of sectionally complemented lattices and the
class of atomistic lattices. Furthermore, it is closed
under direct limit.

\item Let $L$ be a congruence-splitting lattice. Then
$\Conc L$ satisfies \tup{URP}.

\item Every weakly distributive image (in E.~T.
Schmidt's sense, see \cite{Schm68}) of a distributive
lattice with zero satisfies \tup{URP}.

\item $S_\gk$ does not satisfy
\tup{URP}, for all $\gk\geq\aleph_2$.

\item Let $F$ be a free lattice with at least $\aleph_2$
generators in any non-distributive variety of lattices.
Then $\Conc F$ does not satisfy \tup{URP}.
\end{enumerate}
\end{tha}

Observe that every lattice $L$ which is either atomistic
or sectionally complemented has permutable congruences. It
is easy to see that the class of lattices with permutable
congruences is self-dual, and closed under direct limit.

Our following result provides the natural
generalization of Theorem~1, contained in Theorems
\ref{T:CongSplType1}, \ref{T:Type1URP}, \ref{T:URPIMac}, and
Corollary~\ref{C:FreeNonType1.5}:

\begin{thb}
Let $L$ be a lattice.
\begin{enumerate}
\item If $L$ is congruence-splitting, then $L$ has
permutable congruences.

\item Suppose that $L$ has permutable congruences. Then
$\Conc L$ satisfies a certain strengthening of \tup{URP},
denoted here \URPI.

\item Suppose that $L$ has almost permutable congruences,
that is, the equality $\ga\jj\gb=\ga\gb\uu\gb\ga$ holds, for
all congruences $\ga$ and $\gb$ of $L$. Then $\Conc L$ satisfies
a certain weakening of \URPI, denoted here \URPIM.

\item Let $F$ be a free bounded lattice with at least $\aleph_2$
generators in any non-distributive variety of lattices.
Then $\Conc F$ does not satisfy \URPIM.

\end{enumerate}
\end{thb}

On the other hand, the proof of the negative property of the
semilattices $S_\gk$ (or of congruence lattices of free
lattices) given in Theorem 1(iv,v)
relies on a very simple, but powerful, infinite
combinatorial result due to C.~Kuratowski \cite{Kura}, that
characterizes cardinal numbers of the form $\aleph_n$ (and
we just need it for $\aleph_2$) \emph{via} finite-valued set
maps. A look at the proofs suggests that the \emph{finite},
combinatorial core of the problem lies in
\emph{amalgamation properties} of the $\Conc$ functor.

Positive amalgamation results have been obtained by
the first author, who proves, in \cite{Tuma}, that for any
diagram $\E S$ of distributive \jz-semilattices indexed by 
the square $\two^2=\two\times\two$ (with $\two=\set{0,1}$),
any lifting, with respect to the
$\Conc$ functor, of $\E S$ minus
the top semilattice, by finite atomistic lattices can be
extended to a lifting of the full diagram $\E S$, with a
finite atomistic lattice at the top. This somehow cumbersome
formulation is relevant for constructing direct systems of
size $\aleph_1$, see
\cite{GLWe}, where the results of \cite{Tuma} are also
generalized to arbitrary lattices with finite congruence
lattices.

Because of Theorem~1, it is not hard to
see that a $n$-dimensional version of this result, for $n$
arbitrary, cannot hold. Because of the result, stated in
\cite{GS96}, that every finite lattice admits a
congruence-preserving extension into a finite sectionally
complemented lattice, this negative result can be extended
to all finite lattices (not only atomistic).

However, this proof is not constructive, and it remained to
find an effective example of the non-existence of an
amalgamation. This means that one tries to extract a
\emph{combinatorial} information from the existence of
counterexamples of size~$\aleph_2$.

A stronger negative result, also simpler to
state, would be the existence of a diagram of finite
distributive semilattices, indexed, say, by $\two^3$,
without a lifting, with respect to $\Conc$, by finite atomistic
lattices. A finite diagram without a lifting by finite
atomistic lattices has been constructed by the first author
in \cite{Tuma93}, but it is not indexed by any
$n$-dimensional cube. A cube, without amalgamation, of
related objects, called \emph{V-measures}, is constructed in
the paper by H. Dobbertin \cite[pp. 32--34]{Dobb}.

In this paper, we provide a counterexample to the cube
amalgamation problem above:

\goodbreak
\begin{thc}\hfill
\begin{enumerate}
\item There exists a cube of finite Boolean semilattices and
semilattice embeddings, that cannot be lifted,
with respect to the $\Conc$ functor, by lattices with almost
permutable congruences.

\item There exists a cube of finite Boolean semilattices and
semilattice embeddings that can be lifted, with respect to the
$\Conc$ functor, by finite lattices, but that cannot be
lifted, with respect to the $\Conc$ functor, by lattices with
permutable congruences.
\end{enumerate}
\end{thc}

In fact, we provide somewhat stronger statements,
see Theorems \ref{T:NonAmalg}, \ref{T:Lift},
and~\ref{T:NonAmalg2}.

Finally, in Section~\ref{S:NoFunc}, we show a very simple
example of a diagram of finite Boolean semilattices, that
cannot be lifted, in an isomorphism-preserving fashion, with
respect to the $\Conc$ functor, by arbitrary lattices. This
section can be read independently of the others. This shows
that a positive solution of the Congruence Lattice Problem
cannot be achieved by any \emph{functor} from distributive
\jz-semilattices and \jz-homomorphisms, to lattices and
lattice homomorphisms.

\section{Lattices with permutable congruences}\label{S:Comm}

We shall first recall some basic notation and terminology.
If $\ga$ and $\gb$ are two binary relations on a set $A$,
then we put
 \[
 \ga\gb=\setm{\vv<x,y>\in A\times A}
 {(\exists z\in A)(\vv<x,z>\in\ga
 \text{ and }\vv<z,y>\in\gb}.
 \]
If $\ga$ is a binary relation on $A$ and $x$,
$y\in A$, we shall often write $x\equiv_\ga y$ instead of
$\vv<x,y>\in\ga$. Suppose now that $A$ is given a
structure of algebra (over a given similarity type). For
any two elements $x$ and $y$ of $A$, we denote by
$\gQ_A(x,y)$ (or simply $\gQ(x,y)$) the least congruence
of $A$ that identifies $x$ and $y$. Furthermore, we shall
say that $A$ \emph{has permutable congruences}, if for any
congruences $\ga$ and $\gb$ of $A$, the equality
$\ga\gb=\gb\ga$ holds. As usual, the similarity
type of the class of all lattices consists of the two binary
operations $\mm$ and $\jj$.

We shall first give a convenient characterization of
lattices with permutable congruences. The origin of the
argument can be traced back to R.~P. Dilworth's paper
\cite{Dilw}. The authors are grateful to Ralph Freese for
having pointed this out, along with the following proof,
which simplifies tremendously the original one.

\begin{proposition}\label{P:CharType1}
Let $L$ be a lattice. Then the following are equivalent:
\begin{enumerate}
\item $L$ has permutable congruences.

\item For any elements $a$, $b$, and $c$ of $L$ such that
$a\leq c\leq b$, there exists $x\in L$ such that
$a\equiv_{\gQ(c,b)}x$ and $x\equiv_{\gQ(a,c)}b$.
\end{enumerate}

\end{proposition}

Note that if there exists an element $x$ of $L$ as in (ii)
above, then $y=(x\jj a)\mm b$ satisfies that
$a\equiv_{\gQ(c,b)}y$, $y\equiv_{\gQ(a,c)}b$, and
$a\leq y\leq b$. So, the element $x$ of (ii) can be assumed to
lie in the interval $[a,b]$.

\begin{proof}
(i)$\Rightarrow$(ii) Assume that (i) holds. Put
$\ga=\gQ(a,c)$ and $\gb=\gQ(c,b)$. Since $a\equiv_\ga c$
and $c\equiv_\gb b$, we have $\vv<a,b>\in{\ga\gb}$.
By assumption, we also have $\vv<a,b>\in{\gb\ga}$,
which turns out to be the desired conclusion.

(ii)$\Rightarrow$(i) Assume that (ii) holds. Let $\ga$ and
$\gb$ be congruences of $L$, we prove that $\ga\gb$
is contained into $\gb\ga$. Thus let $a$, $b$,
$c\in L$ such that $a\equiv_{\ga}b$ and $b\equiv_{\gb}c$.
Then $a\equiv_{\ga}a\jj b$ and
$a\jj b\equiv_{\gb}a\jj b\jj c$, thus, by the assumption,
there exists $x$ such that
$a\equiv_{\gQ(a\jj b,a\jj b\jj c)}x$ and
$x\equiv_{\gQ(a,a\jj b)}a\jj b\jj c$. Thus
$a\equiv_{\gb}x$ and
$x\equiv_{\ga}a\jj b\jj c$. Similarly, by reversing the
roles of $a$ and $c$ as well as $\ga$ and $\gb$, there
exists $y$ such that $c\equiv_{\ga}y$ and
$y\equiv_{\gb}a\jj b\jj c$. Put $z=x\mm y$. Then
$a\equiv_{\gb}z$ and $z\equiv_{\ga}c$; whence
$a\equiv_{\gb\ga}c$.
\end{proof}

We recall now the following definition of \cite{Wehr2},
also used in \cite{PTWe}. A lattice $L$ is
\emph{congruence-splitting}, if for all $u\leq v$ in $L$
and all $\ga$, $\gb\in\Con L$ such that
$\gQ(u,v)=\ga\jj\gb$, there exist $x$ and $y$ in the
interval $[u,v]$ such that $x\jj y=v$, $u\equiv_\ga x$,
and $u\equiv_\gb y$.

As an immediate corollary of
Proposition~\ref{P:CharType1}, we obtain the following:

\begin{theorem}\label{T:CongSplType1}
Every congruence-splitting lattice has permutable
congruences.
\end{theorem}

\begin{proof}
Let $L$ be a congruence-splitting lattice. We shall prove
that the assumption (ii) of Proposition~\ref{P:CharType1}
holds. Thus let $a\leq c\leq b$ in $L$. Put $\ga=\gQ(a,c)$
and $\gb=\gQ(c,b)$. Since $\ga\jj\gb=\gQ(a,b)$, there
exist, by assumption, elements $x$ and $y$ of $[a,b]$ such
that $x\jj y=b$, and $a\equiv_\ga x$ and $a\equiv_\gb y$.
Joining with $y$ the first relation, we obtain
that $y\equiv_\ga b$. So, $a\equiv_\gb y$ and
$y\equiv_\ga b$, which is the desired conclusion.
\end{proof}

Congruence-splitting lattices have been introduced in
\cite{Wehr2} because their compact congruence semilattices
satisfy a certain infinite axiom, called the \emph{Uniform
Refinement Property}, denoted URP. We introduce here an
apparently stronger property, denoted \URPI, and we prove
that the compact congruence semilattice of every lattice
with permutable congruences satisfies this property.

\begin{definition}\label{D:URPI}
Let $S$ be a \jz-semilattice, and let $\ge$ be an element
of $S$. We say that $S$ \emph{satisfies \URPI\ at $\ge$},
if for every family $\vv<\vv<\ga_i,\gb_i>\mid i\in I>$ of
elements of $S\times S$ such that $\ga_i\jj\gb_i=\ge$ for
all $i$, there exist a family
$\vv<\vv<\ga_i^*,\gb_i^*>\mid i\in I>$ of elements of
$S\times S$ and a family
$\vv<\gg_{i,j}\mid\vv<i,j>\in I\times I>$
of elements of $S$ such that for
all $i$, $j$, $k\in I$, the following conditions hold:
\begin{enumerate}
\item $\ga_i^*\leq\ga_i$ and $\gb_i^*\leq\gb_i$, and
$\ga_i^*\jj\gb_i^*=\ge$.

\item $\gg_{i,j}\leq\ga_i^*$ and $\gg_{i,j}\leq\gb_j^*$.

\item $\ga_i^*\leq\ga_j^*\jj\gg_{i,j}$ and
$\gb_j^*\leq\gb_i^*\jj\gg_{i,j}$.

\item $\gg_{i,k}\leq\gg_{i,j}\jj\gg_{j,k}$.
\end{enumerate}
Furthermore, say that $S$ satisfies \URPI, if $S$
satisfies \URPI\ at every element of $S$.
\end{definition}

It may appear strange at first sight to refer to the
system of the $\ga_i^*$, $\gb_i^*$, $\gg_{i,j}$ as a
\emph{refinement} of the $\ga_i$, $\gb_i$. Section~2
of~\cite{Wehr2} gives some motivation for this terminology.

The proof of the following result is essentially the same
as the proof of Proposition~2.2 of \cite{Wehr2}, so we
shall omit it:

\begin{proposition}\label{P:URPAdd}
Let $S$ be a distributive semilattice. Then the set of
all elements of $S$ at which \URPI\ holds is
closed under the join operation.
\end{proposition}

\begin{theorem}\label{T:Type1URP}
Let $L$ be a lattice with permutable congruences.
Then $\Conc L$ satisfies \URPI.
\end{theorem}

\begin{proof}
Since $\Conc L$ is a distributive semilattice, it suffices,
by Proposition~\ref{P:URPAdd}, to prove that $\Conc L$
satisfies \URPI\ at every element of the form
$\ge=\gQ(u,v)$, where $u\leq v$ in $L$. Thus consider a
family $\vv<\vv<\ga_i,\gb_i>\mid i\in I>$ of elements of
$\Conc L$ such that $\ga_i\jj\gb_i=\ge$, for all $i\in I$.
Since $L$ has permutable congruences, we have also
$\vv<u,v>\in\ga_i\gb_i$, for all $i\in I$, thus there
exists $x_i\in[u,v]$ such that $u\equiv_{\ga_i}x_i$ and
$x_i\equiv_{\gb_i}v$. By replacing $x_i$ by
$(u\jj x_i)\mm v$, one can suppose, without loss of
generality, that $u\leq x_i\leq v$. Put
$\ga_i^*=\gQ(u,x_i)$, $\gb_i^*=\gQ(x_i,v)$ and
$\gg_{i,j}=\gQ(x_i,x_i\mm x_j)$, for all $i$, $j\in I$. It
is easy to verify that the elements $\ga_i^*$, $\gb_i^*$,
$\gg_{i,j}$ verify the relations (i) to (iv) of the
definition of \URPI\ (Definition~\ref{D:URPI}).
\end{proof}

It is proved in Corollary~4.1 of \cite{PTWe} that if
$\VV$ is any non-distributive variety of lattices
and if $F$ is a free lattice in $\VV$ on at least
$\aleph_2$ generators, then $\Conc F$ does not satisfy a
certain weakening of \URPI, denoted there WURP. As a
corollary, we note, in particular, the following:

\begin{corollary}\label{C:FreeNonType1}
Let $\VV$ be any non-distributive variety of
lattices, and let $F$ be any free lattice in $\VV$
on at least $\aleph_2$ generators. Then there is no
lattice $K$ with permutable congruences such that
$\Conc K\iso\Conc F$.
\end{corollary}

\section{Lattices with almost permutable congruences}
\label{S:AComm}

Two congruences $\ga$ and $\gb$ of an algebra $A$ are said to
be \emph{almost permutable}, if the following equality holds:
 \[
 \ga\jj\gb=\ga\gb\uu\gb\ga.
 \]
We say that $A$ has \emph{almost permutable congruences}, if
any two congruences of $A$ are almost permutable. The
three-element chain is an easy example of a lattice with almost
permutable congruences but not with permutable congruences.

We shall formulate an analogue of Definition~\ref{D:URPI} for
lattices with almost permutable congruences:

\begin{definition}\label{D:URPIM}
Let $S$ be a \jz-semilattice, let $\ge\in S$, and let
 \[
 \gs=\vv<\vv<\ga_i,\gb_i>\mid i\in I>
 \]
be a family of elements of $S\times S$ such that
 \begin{equation}\label{Eq:gagbge}
 \ga_i\jj\gb_i=\ge,\text{ for all }i\in I.
 \end{equation}
We say that $S$ \emph{satisfies \URPIM\ at $\gs$},
if there exist a subset $X$ of $I$, a family
$\vv<\vv<\ga_i^*,\gb_i^*>\mid i\in I>$ of elements of
$S\times S$ and a family
$\vv<\gg_{i,j}\mid\vv<i,j>\in I\times I>$
of elements of $S$ such that for
all $i$, $j$, $k\in I$, the following conditions hold:
\begin{enumerate}
\item $\ga_i^*\leq\ga_i$, $\gb_i^*\leq\gb_i$, and
$\ga_i^*\jj\gb_i^*=\ge$.

\item $\gg_{i,j}\leq\ga_i^*$ and $\gg_{i,j}\leq\gb_j^*$.

\item $\ga_i^*\leq\ga_j^*\jj\gg_{i,j}$ and
$\gb_j^*\leq\gb_i^*\jj\gg_{i,j}$.

\item $\gg_{i,k}\leq\gg_{i,j}\jj\gg_{j,k}$,
for all $i$, $j$, $k\in I$ such that
 \[
 \set{i,k}\ci X\Rightarrow j\in X\text{ and }
 \set{i,k}\ci I\setminus X\Rightarrow j\in I\setminus X.
 \]
\end{enumerate}
We say that $S$ satisfies \URPIM\ at an element $\ge$ of $S$,
if $S$ satisfies \URPIM\ at every family $\gs$ satisfying
\eqref{Eq:gagbge}. Finally, we say that $S$ satisfies \URPIM,
if it satisfies \URPIM\ at $\ge$, for all $\ge\in S$.
\end{definition}

We recall that if $S$ and $T$ are join-semilattices and if
$\ge\in S$, a \jhom\ $\gm\colon S\to T$ is
\emph{weakly distributive} at $\ge$, if for all $\ga$,
$\gb\in T$ such that $\ga\jj\gb=\gm(\ge)$, there are $\ga'$,
$\gb'\in S$ such that the following relations hold:
 \[
 \ge=\ga'\jj\gb',\quad\gm(\ga')\leq\ga,\quad\gm(\gb')\leq\gb.
 \]
We now formulate an analogue of Proposition~2.3 of
\cite{Wehr2}. The proof is straightforward, thus we omit it.

\begin{lemma}\label{L:wdimage}
Let $S$ and $T$ be join-semilattices, let $\ge\in S$.
Let $\gm\colon S\to T$ be a weakly
distributive \jhom. If $S$ satisfies \URPIM\ at
$\ge$, then $T$ satisfies \URPIM\ at $\gm(\ge)$.
\end{lemma}

We can now use the results of \cite{PTWe} to get the following
negative result:

\begin{corollary}\label{C:FreeNonType1.5}
Let $\VV$ be any non-distributive variety of
lattices, and let $F$ be any free bounded lattice in $\VV$
on at least $\aleph_2$ generators. Then $\Conc F$ does not
satisfy \URPIM\ at the largest congruence, $\gQ(0,1)$,
of $F$.
\end{corollary}

\begin{proof}
As in \cite{PTWe}, we shall denote by $\BB_\VV(I)$ the free
product (=coproduct) of $I$ copies of the two-element chain
in $\VV$, say, $s_i<t_i$, for $i\in I$, with bounds $0$, $1$
added. Now assume that $|I|=\aleph_2$. As in Corollary~4.1 in
\cite{PTWe}, $\BB_\VV(I)$ is a quotient of $F$, thus,
by Proposition~1.2 of \cite{Wehr2}, the induced map $f$ from
$\Conc F$ onto
$\Conc\BB_\VV(I)$ is weakly distributive. Since
$f(\gQ_F(0,1))$ equals $\gQ_{\BB_\VV(I)}(0,1)$, it suffices,
by Lemma~\ref{L:wdimage}, to prove that $\Conc\BB_\VV(I)$ does
not satisfy \URPIM\ at $\gQ_{\BB_\VV(I)}(0,1)$.

We follow for this the pattern of the proof that begins in
Section~2 of \cite{PTWe}. Suppose that $\Conc\BB_\VV(I)$
satisfies \URPIM\ at $\gQ_{\BB_\VV(I)}(0,1)$. We put, again,
 \begin{gather}
 \ga_i=\gQ(0,s_i)\jj\gQ(t_i,1)\quad\text{and}\quad
 \gb_i=\gQ(s_i,t_i),
 \quad\text{for all }i\in I,\label{Eq:Defgaigbi}\\
 \ge=\gQ(0,1).\label{Eq:Defeps}
 \end{gather}
Observe that $\ga_i\jj\gb_i=\ge$, for all $i\in I$. By using
\URPIM, we obtain a subset $X$ of $I$ and elements $\gg_{i,j}$ of
$\Conc\BB_\VV(I)$, for $i$,
$j\in I$, satisfying the following conditions:
 \begin{align}
 \gg_{i,j}\ci\ga_i,\,\gb_j,&
 \quad\text{for all }i,\,j\in I;\label{Eq:ref1}\\
 \gg_{i,j}\jj\ga_j\jj\gb_i=\ge,&
 \quad\text{for all }i,\,j\in I;\label{Eq:ref2}\\
 \gg_{i,k}\ci\gg_{i,j}\jj\gg_{j,k},&
 \quad\text{for all }i,\,j,\,k\in I\text{ such that}\notag\\
 &\quad\set{i,k}\ci X\Rightarrow j\in X\text{ and }
 \set{i,k}\ci I\setminus X\Rightarrow j\in I\setminus X.
 \label{Eq:ref3}
 \end{align}
Note that either $|X|=\aleph_2$ or $|I\setminus X|=\aleph_2$,
say, without loss of generality, $|X|=\aleph_2$. By restricting
\eqref{Eq:ref1}--\eqref{Eq:ref3} to $X$, we obtain the
following conditions:
 \begin{align*}
 \gg_{i,j}\ci\ga_i,\,\gb_j,&
 \quad\text{for all }i,\,j\in X;\\
 \gg_{i,j}\jj\ga_j\jj\gb_i=\ge,&
 \quad\text{for all }i,\,j\in X;\\
 \gg_{i,k}\ci\gg_{i,j}\jj\gg_{j,k},&
 \quad\text{for all }i,\,j,\,k\in X.
 \end{align*}
These relations hold in the join-semilattice
$\Conc\BB_\VV(I)$, however, by projecting them onto
$\Conc\BB_\VV(X)$ (use any lattice retraction from $\BB_\VV(I)$ 
onto $\BB_\VV(X)$ that preserves the bounds and the $s_i$, $t_i$,
for $i\in X$), we obtain relations of the form
 \begin{align*}
 \gg'_{i,j}\ci\ga_i,\,\gb_j,&
 \quad\text{for all }i,\,j\in X;\\
 \gg'_{i,j}\jj\ga_j\jj\gb_i=\ge,&
 \quad\text{for all }i,\,j\in X;\\
 \gg'_{i,k}\ci\gg'_{i,j}\jj\gg'_{j,k},&
 \quad\text{for all }i,\,j,\,k\in X,
 \end{align*}
for elements $\gg'_{i,j}$ of $\Conc\BB_\VV(X)$, for
$i$, $j\in X$. In the formulas above, we slightly abuse the
notation by still denoting $\ga_i$, $\gb_i$, and $\ge$ the compact
congruences of $\BB_\VV(X)$ (not $\BB_\VV(I)$) defined within
$\BB_\VV(X)$ by the formulas \eqref{Eq:Defgaigbi} and
\eqref{Eq:Defeps}.
However, since $|X|=\aleph_2$, this cannot exist by the
proof of Theorem~3.3 of \cite{PTWe}.
\end{proof}

On the positive side, the following result relates \URPIM\
with lattices with almost permutable congruences:

\begin{theorem}\label{T:URPIMac}
Let $L$ be a lattice. If $L$ has almost permutable congruences,
then $\Conc L$ satisfies \URPIM\ at every principal congruence
of $L$.
\end{theorem}

\begin{proof}
Let $u\leq v$ in $L$, we prove that $\Conc L$ satisfies
\URPIM\ at $\ge=\gQ_L(u,v)$. So, let
$\vv<\vv<\ga_i,\gb_i>\mid i\in I>$ be a family of ordered pairs
of compact congruences of $L$ such that $\ga_i\jj\gb_i=\ge$,
for all $i\in I$. Since $L$ has almost permutable congruences,
there exists a subset $X$ of $I$ such that
 \begin{align*}
 u&\equiv_{\ga_i\gb_i}v,
 \quad\text{for all }i\in X,\\
 u&\equiv_{\gb_i\ga_i}v,
 \quad\text{for all }i\in I\setminus X.
 \end{align*}
Thus, for all $i\in I$, we obtain an element $x_i$ of the
interval $[u,v]$ such that
 \begin{alignat*}{2}
 u&\equiv_{\ga_i}x_i&\text{ and }x_i\equiv_{\gb_i}v,
 &\quad\text{for all }i\in X,\\
 u&\equiv_{\gb_i}x_i&\text{ and }x_i\equiv_{\ga_i}v,
 &\quad\text{for all }i\in I\setminus X.
 \end{alignat*}
We define compact congruences $\ga_i^*$ and $\gb_i^*$
of $L$, for all $i\in I$, as follows:
 \begin{alignat*}{2}
 \ga_i^*&=\gQ(u,x_i)&\text{ and }\gb_i^*=\gQ(x_i,v),&\quad
 \text{for all }i\in X,\\
 \ga_i^*&=\gQ(x_i,v)&\text{ and }\gb_i^*=\gQ(u,x_i),&\quad
 \text{for all }i\in I\setminus X.
 \end{alignat*}
Note the following properties of $\ga_i^*$ and $\gb_i^*$:
\begin{gather}
\ga_i^*\ci\ga_i;\quad\gb_i^*\ci\gb_i;\\
\ga_i^*\jj\gb_i^*=\ge,
\end{gather}
for all $i\in I$.
We further define compact congruences $\gg_{i,j}$ of $L$, for
$i$, $j\in I$, as follows:
 \begin{equation*}
 \gg_{i,j}=
  \begin{cases}
  \gQ^+(x_i,x_j),&\text{if }i\in X\text{ and }j\in X;\\
  \gQ(u,x_i\mm x_j),&\text{if }i\in X\text{ and }j\nin X;\\
  \gQ(x_i\jj x_j,v),&\text{if }i\nin X\text{ and }j\in X;\\
  \gQ^+(x_j,x_i),&\text{if }i\nin X\text{ and }j\nin X.\\
  \end{cases} 
 \end{equation*}
We use here the convenient notation
 \[
 \gQ^+(a,b)=\gQ(a\mm b,a),\text{ for all }a,\,b\in L.
 \]
Next, we verify that the congruences $\ga_i^*$, $\gb_i^*$, and
$\gg_{i,j}$ satisfy the relations listed in
Definition~\ref{D:URPIM}. We start with the following
 \begin{equation}\label{Eq:ciaibig}
 \gg_{i,j}\ci\ga_i^*,\,\gb_j^*,
 \quad\text{for all }i,\,j\in I.
 \end{equation}
There are four cases to consider. If $i$, $j\in X$, then the
verification of \eqref{Eq:ciaibig} amounts to the verification
of the containment
 \[
 \gQ^+(x_i,x_j)\ci\gQ(u,x_i),\,\gQ(x_j,v),
 \]
which is immediate. Similarly, the remaining cases $i\in X$
and $j\nin X$, $i\nin X$ and $j\in X$, and $i$, $j\nin X$,
respectively, correspond to the containments
 \begin{align*}
 \gQ(u,x_i\mm x_j)&\ci\gQ(u,x_i),\,\gQ(u,x_j),\\
 \gQ(x_i\jj x_j,v)&\ci\gQ(x_i,v),\,\gQ(x_j,v),\\
 \gQ^+(x_j,x_i)&\ci\gQ(x_i,v),\,\gQ(u,x_j), 
 \end{align*}
all of which are obvious. This completes the verification of
\eqref{Eq:ciaibig}.

We proceed with the verification of
 \begin{equation}\label{Eq:ceaibig1}
 \ga_i^*\ci\ga_j^*\jj\gg_{i,j}
 \quad\text{for all }i,\,j\in I.
 \end{equation}
Again, there are the same four cases as before to consider,
that correspond respectively to the following containments:
 \begin{align*}
 \gQ(u,x_i)&\ci\gQ(u,x_j)\jj\gQ^+(x_i,x_j),\\
 \gQ(u,x_i)&\ci\gQ(x_j,v)\jj\gQ(u,x_i\mm x_j),\\
 \gQ(x_i,v)&\ci\gQ(u,x_j)\jj\gQ(x_i\jj x_j,v),\\
 \gQ(x_i,v)&\ci\gQ(x_j,v)\jj\gQ^+(x_j,x_i),
 \end{align*}
all of which are obvious, thus completing the verification
of \eqref{Eq:ceaibig1}.

The following system of containments can be verified in a similar
fashion:
 \begin{equation}\label{Eq:ceaibig2}
 \gb_j^*\ci\gb_i^*\jj\gg_{i,j}
 \quad\text{for all }i,\,j\in I.
 \end{equation}

Finally, let $i$, $j$, $k\in I$, and
suppose that $\set{i,k}\ci X$ implies $j\in X$, and
$\set{i,k}\ci I\setminus X$ implies $j\in I\setminus X$. We
prove that $\gg_{i,k}\ci\gg_{i,j}\jj\gg_{j,k}$. There are six
cases to consider, which are, respectively,
 \begin{align*}
 i\in X;\ j\in X;\ k\in X,\\
 i\in X;\ j\in X;\ k\nin X,\\
 i\in X;\ j\nin X;\ k\nin X,\\
 i\nin X;\ j\in X;\ k\in X,\\
 i\nin X;\ j\nin X;\ k\in X,\\
 i\nin X;\ j\nin X;\ k\nin X.
 \end{align*}
The corresponding containments to be verified are
 \begin{align*}
 \gQ^+(x_i,x_k)&\ci\gQ^+(x_i,x_j)\jj\gQ^+(x_j,x_k),\\
 \gQ(u,x_i\mm x_k)&\ci\gQ^+(x_i,x_j)\jj\gQ(u,x_j\mm x_k),\\
 \gQ(u,x_i\mm x_k)&\ci\gQ(u,x_i\mm x_j)\jj\gQ^+(x_k,x_j),\\
 \gQ(x_i\jj x_k,v)&\ci\gQ(x_i\jj x_j,v)\jj\gQ^+(x_j,x_k),\\
 \gQ(x_i\jj x_k,v)&\ci\gQ^+(x_j,x_i)\jj\gQ(x_j\jj x_k,v),\\
 \gQ^+(x_k,x_i)&\ci\gQ^+(x_j,x_i)\jj\gQ^+(x_k,x_j),
 \end{align*}
all of which are obvious, thus completing the proof.
\end{proof}

By putting together the results of Theorem~\ref{T:URPIMac} and
Corollary~\ref{C:FreeNonType1.5}, we thus obtain the following
result:

\begin{corollary}\label{C:NonReprAC}
Let $\VV$ be a non-distributive variety of lattices, and let
$F$ be any free bounded lattice in $\VV$ on at least
$\aleph_2$ generators. Then there exists no lattice $L$ with
almost permutable congruences such that $\Conc L\iso\Conc F$.
\end{corollary}

\begin{remark}
We could have simplified the statement of \URPIM\ in
Definition~\ref{D:URPIM}, and thus the proof of
Theorem~\ref{T:URPIMac}, to obtain exactly exactly the same
statement of Corollary~\ref{C:NonReprAC}. For example, the
statements (i)--(iv) of Definition~\ref{D:URPIM} could
have been simplified the same way as in \cite{PTWe},
namely, into
\begin{itemize}
\item[(ii$'$)] $\gg_{i,j}\ci\ga_i$, $\gb_j$,
for all $i$, $j\in I$,

\item[(iii$'$)] $\gg_{i,j}\jj\ga_j\jj\gb_i=\ge$,
for all $i$, $j\in I$,

\item[(iv$'$)] $\gg_{i,k}\ci\gg_{i,j}\jj\gg_{j,k}$,
for all $i$, $j$, $k\in I$ such that
 \[
 \text{either }\set{i,j,k}\ci X\text{ or }
 \set{i,j,k}\ci I\setminus X.
 \]
\end{itemize}

However, this would have made us miss an important point:
namely, that part (iv) of the original definition of
\URPIM\ ``almost'' holds \emph{for all} $i$, $j$, $k\in X$.
The offending cases are, respectively,
 \begin{align*}
 i\in X;\ j\nin X;\ k\in X,\\
 i\nin X;\ j\in X;\ k\nin X, 
 \end{align*}
and the corresponding containments are, respectively,
 \begin{align*}
 \gQ^+(x_i,x_k)&\ci\gQ(0,x_i\mm x_j)\jj\gQ(x_j\jj x_k,1),\\
 \gQ^+(x_k,x_i)&\ci\gQ(x_i\jj x_j,1)\jj\gQ(x_j\jj x_k,1), 
 \end{align*}
which are easily seen to fail in very simple finite lattices,
for example, the five-element modular non-distributive lattice
$M_3$, by assigning to $x_i$, $x_j$, and $x_k$ the
three atoms of $M_3$. This suggests that \URPI\ may, in fact,
\emph{not} hold in general for the semilattices $\Conc L$,
where $L$ has almost permutable congruences. This suggests, in
the longer term, that no ``uniform refinement property'' of any
sort holds for the semilattices $\Conc L$, where $L$ is an
arbitrary lattice.
\end{remark}

\section{The basic semilattice diagram, $\Dc$}\label{S:SemCube}

We shall construct in this section a finite diagram,
$\Dc$, of finite distributive \jz-semilattices. The ultimate
purpose of this construction will be completed in
Section~\ref{S:NoLift}, where it will be proved that $\Dc$ has
no lifting, with respect to the $\Conc$ functor, by lattices
with permutable congruences. The semilattices of $\Dc$ are
\emph{free} semilattices. If $S$ is a
\jz-semilattice and if $n\in\go$, a \emph{free $n$-tuple} of
elements of $S$ is an element $\vv<s_i\mid i<n>$ of $S^n$
such that the map from the powerset semilattice
$\vv<\E P(n),\jj,\es>$ to $S$ defined by the rule
 \[
 X\mapsto\JJm{s_i}{i\in X}
 \]
is an \emph{embedding} of \jz-semilattices.

For any poset $\vv<P,\leq>$, we view $P$ as a category in the
usual fashion, that is, the objects are the elements of $P$,
and, for $p$, $q\in P$, there exists exactly one morphism
from $p$ to $q$ if $p\leq q$, and none otherwise.
Let $\AA$ be a category. A \emph{$P$-diagram} of $\AA$ is
a functor from $P$ (or, more precisely, the
category associated with $P$) to $\AA$.

Now let $\AA$ and $\BB$ be two categories, let $\FF$ be a
functor from $\AA$ to $\BB$, and let $P$ be a poset. A
$P$-diagram $f$ on $\AA$ is a \emph{lifting} of a
$P$-diagram $g$ on $\BB$, if both diagrams $\FF\circ f$
and $g$ are isomorphic, in notation, $\FF\circ f\iso g$.

The case in which we shall be interested is here the
following. The category $\AA$ is the category of all
lattices and lattice homomorphisms, the category $\BB$ is the
category of all distributive \jz-semilattices and
\jz-homomorphisms, and $\FF$ is the $\Conc$ functor from $\AA$
to $\BB$.
\smallskip

We shall now define \jz-semilattices $S_0$, $S_1$, $S_2$,
$T_0$, $T_1$, $T_2$, and $U$, the building stones of $\Dc$.

\begin{enumerate}
\item $U$ is the powerset semilattice of the five-element
set $5=\set{0,1,2,3,4}$, $U=\mathcal{P}(5)$.

We define elements $\gx_i$, $\gh_i$, and $\gz_i$ ($i<4$) of
$U$ as follows:
 \begin{align*}
 \gx_0&=\set{0,4},&\gx_1&=\set{3},&
 \gx_2&=\set{2},&\gx_3&=\set{1,4};\\
 \gh_0&=\set{0,4},&\gh_1&=\set{1,4},&
 \gh_2&=\set{2},&\gh_3&=\set{3,4};\\
 \gz_0&=\set{0,4},&\gz_1&=\set{1},&
 \gz_2&=\set{3},&\gz_3&=\set{2,4}.
 \end{align*}
Furthermore, we shall denote by $\zero$ (resp., $\one$) the
smallest (resp., largest) element of $U$. Note that
the following equalities hold:
 \[
 \one=\JJm{\gx_i}{i<4}=\JJm{\gh_i}{i<4}=\JJm{\gz_i}{i<4}.
 \]

We shall now define certain subsemilattices of $U$.

\item We define $T_0$ to be the \jz-subsemilattice of $U$
generated by $\setm{\gx_j}{j<4}$. Similarly,
define $T_1$ to be the \jz-subsemilattice of $U$
generated by\linebreak
$\setm{\gh_j}{j<4}$, and $T_2$ to be the \jz-subsemilattice
of $U$ generated by $\setm{\gz_j}{j<4}$.

\item Finally, for all $i<3$, let $S_i$ be the
\jz-subsemilattice of $U$ generated by $\set{\ga_i,\gb_i}$,
where we put
 \begin{align*}
 \ga_0&=\set{0,1,4},&\gb_0&=\set{2,3,4};\\
 \ga_1&=\set{0,3,4},&\gb_1&=\set{1,2,4};\\
 \ga_2&=\set{0,2,4},&\gb_2&=\set{1,3,4}.
 \end{align*}
\end{enumerate}

Note, in particular, that $\one$ is the largest element of
$S_i$ and of $T_i$, for all $i<3$. Furthermore, note that
the definitions of the $\ga_i$ and $\gb_i$ stated in (iii)
imply immediately that the following arrays

 \begin{equation}\label{Eq:RefMat}
 \mbox{\begin{tabular}{|c|c|c|}
 \cline{2-3}
 \multicolumn{1}{l|}{} & $\tvi \ga_2$ & $\gb_2$\\
 \hline
 $\tvi \ga_1$ & $\gx_0$ & $\gx_1$\\
 \hline
 $\tvi \gb_1$ & $\gx_2$ & $\gx_3$\\
 \hline
 \end{tabular}}\qquad
 \mbox{\begin{tabular}{|c|c|c|}
 \cline{2-3}
 \multicolumn{1}{l|}{} & \tvi $\ga_2$ & $\gb_2$\\
 \hline
 $\tvi \ga_0$ & $\gh_0$ & $\gh_1$\\
 \hline
 $\tvi \gb_0$ & $\gh_2$ & $\gh_3$\\
 \hline
 \end{tabular}}\qquad
 \mbox{\begin{tabular}{|c|c|c|}
 \cline{2-3}
 \multicolumn{1}{l|}{} & \tvi $\ga_1$ & $\gb_1$\\
 \hline
 \tvi $\ga_0$ & $\gz_0$ & $\gz_1$\\
 \hline
 \tvi $\gb_0$ & $\gz_2$ & $\gz_3$\\
 \hline
 \end{tabular}}
 \end{equation}
are \emph{refinement matrices}, that is, in each of the
three arrays, the first element of each row is the join of the
other two, and similarly for the columns. For example,
$\ga_0=\gz_0\jj\gz_1=\gh_0\jj\gh_1$,
$\gb_1=\gx_2\jj\gx_3=\gz_1\jj\gz_3$, \emph{etc.}.

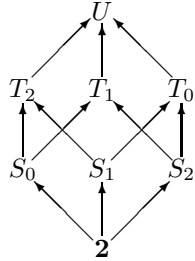
\begin{figure}[htb]
\begin{picture}(100,90)(0,0)

\put(30,0){\makebox(0,0){$\two$}}
\put(0,30){\makebox(0,0){$S_0$}}
\put(30,30){\makebox(0,0){$S_1$}}
\put(60,30){\makebox(0,0){$S_2$}}
\put(0,60){\makebox(0,0){$T_2$}}
\put(30,60){\makebox(0,0){$T_1$}}
\put(60,60){\makebox(0,0){$T_0$}}
\put(30,90){\makebox(0,0){$U$}}

\put(30,5){\vector(0,1){20}}
\put(26.5,3.5){\vector(-1,1){22}}
\put(33.5,3.5){\vector(1,1){22}}
\put(0,35){\vector(0,1){20}}
\put(3.5,33.5){\vector(1,1){22}}
\put(26.5,33.5){\vector(-1,1){22}}
\put(33.5,33.5){\vector(1,1){22}}
\put(60,35){\vector(0,1){20}}
\put(56.5,33.5){\vector(-1,1){22}}
\put(3.5,65){\vector(1,1){22}}
\put(30,65){\vector(0,1){20}}
\put(56.5,65){\vector(-1,1){22}}

\end{picture}
\caption{Semilattice diagram}\label{Fi:DiagSemil}
\end{figure}

At the bottom of the construction, we put the two element
\jz-semilattice, $\two=\vv<\set{0,1},\jj,0>$.
We can see right away that if $i\ne j$ are elements of
$3=\set{0,1,2}$, then
$\two\subseteq S_i\subseteq T_j\subseteq U$. Thus the
semilattices $\two$, $S_i$, $T_i$ ($i<3$), and $U$ can be
arranged in a commutative diagram, as shows
Figure~\ref{Fi:DiagSemil}, where the arrows represent
the inclusion maps. Note that all the maps in
Figure~\ref{Fi:DiagSemil} are embeddings of
\jz-semilattices.

The proof of the following lemma is trivial.

\begin{lemma}\label{L:FreeBases}
The quadruples $\vv<\gx_0,\gx_1,\gx_2,\gx_3>$ (resp.,
$\vv<\gh_0,\gh_1,\gh_2,\gh_3>$,
$\vv<\gz_0,\gz_1,\gz_2,\gz_3>$) are free. Therefore, the
following isomorphisms hold:
\begin{enumerate}
\item $T_0\iso T_1\iso T_2\iso\two^4$.

\item $S_0\iso S_1\iso S_2\iso\two^2$.
\end{enumerate}
In particular, the $S_i$ and the $T_i$, $i<3$, are
finite Boolean semilattices.
\end{lemma}

\begin{lemma}\label{L:NonIneq}
$\gh_1\nleq\gx_1\jj\gz_1$.
\end{lemma}

\begin{proof}
Note that $4\in\gh_1$, but that $4\notin\gx_1\jj\gz_1$.
\end{proof}

\section{Non-existence of a lifting with permutable congruences}
\label{S:NoLift}

We shall prove in this section our first negative lifting
result:

\begin{theorem}\label{T:NonAmalg}
There is no lifting, with respect to the $\Conc$ functor, in
the category of lattices, of the diagram $\Dc$, such that the
lattices corresponding to $S_i$, $i<3$, have permutable
congruences.
\end{theorem}

\begin{proof}
Suppose otherwise. Let us consider a lifting, with respect to
the $\Conc$ functor, of $\Dc$,
by a lattice diagram as in Figure~\ref{Fi:DiagLatt}.

\begin{figure}[hbt]
\begin{picture}(100,140)(0,10)

\put(50,0){\makebox(0,0){$K$}}
\put(0,50){\makebox(0,0){$K_0$}}
\put(50,50){\makebox(0,0){$K_1$}}
\put(100,50){\makebox(0,0){$K_2$}}
\put(0,100){\makebox(0,0){$L_2$}}
\put(50,100){\makebox(0,0){$L_1$}}
\put(100,100){\makebox(0,0){$L_0$}}
\put(50,150){\makebox(0,0){$P$}}

\put(15,30){\makebox(0,0)[r]{$f_0$}}
\put(45,30){\makebox(0,0)[r]{$f_1$}}
\put(85,30){\makebox(0,0)[l]{$f_2$}}

\put(-2,80){\makebox(0,0)[r]{$g_{02}$}}
\put(46,80){\makebox(0,0)[r]{$g_{01}$}}

\put(26,88){\makebox(0,0)[br]{$g_{12}$}}
\put(75,88){\makebox(0,0)[bl]{$g_{10}$}}

\put(52,80){\makebox(0,0)[l]{$g_{21}$}}
\put(102,80){\makebox(0,0)[l]{$g_{20}$}}

\put(75,126){\makebox(0,0)[bl]{$h_0$}}
\put(48,126){\makebox(0,0)[br]{$h_1$}}
\put(24,126){\makebox(0,0)[br]{$h_2$}}

\put(50,5){\vector(0,1){40}}
\put(46.5,3.5){\vector(-1,1){42}}
\put(53.5,3.5){\vector(1,1){43}}
\put(0,55){\vector(0,1){40}}
\put(3.5,53.5){\vector(1,1){42}}
\put(46.5,53.5){\vector(-1,1){42}}
\put(53.5,53.5){\vector(1,1){42}}
\put(100,55){\vector(0,1){40}}
\put(96.5,53.5){\vector(-1,1){42}}
\put(3.5,105){\vector(1,1){42}}
\put(50,105){\vector(0,1){40}}
\put(96.5,105){\vector(-1,1){42}}

\end{picture}
\caption{Lattice diagram}\label{Fi:DiagLatt}
\end{figure}
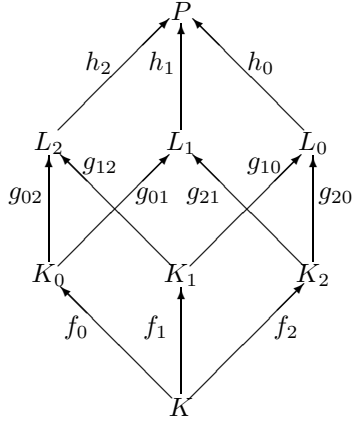

Moreover, suppose that $K_0$, $K_1$, and $K_2$ have
permutable congruences. We shall obtain a contradiction.

As on Figure~\ref{Fi:DiagLatt}, denote by
$f_i\colon K\to K_i$,
$g_{ij}\colon K_i\to L_j$ and $h_j\colon L_j\to P$ the
lattice homomorphisms from the lattice diagram of
Figure~\ref{Fi:DiagLatt}, for all $i$, $j<3$ such that
$i\ne j$. Furthermore, let
$g_i\colon K_i\to P$ be the homomorphism defined by
$g_i=h_j\circ g_{ij}$, for all $j<3$ such that $i\ne j$.
For each lattice $L$ of Figure~\ref{Fi:DiagLatt}, let
$\gm_L$ be the isomorphism from $\Conc L$ onto the
corresponding semilattice of Figure~\ref{Fi:DiagSemil},
such that the isomorphisms $\gm_K\colon\Conc K\to\two$,
$\gm_{K_i}\colon\Conc K_i\to S_i$ (for $i<3$),
$\gm_{L_j}\colon\Conc L_j\to T_j$ (for $j<3$),
$\gm_P\colon\Conc P\to U$ witness the isomorphism of $\Dc$ and
the image by the
$\Conc$ functor of the diagram of Figure~\ref{Fi:DiagLatt}.
For each of those lattices $L$,
put $\gF_L(x,y)=\gm_L\bigl(\gQ_L(x,y)\bigr)$ (so, an element
of $\gm_L\bigl[\Conc L\bigr]$), for all $x$, $y\in L$. For
example, if $x$, $y\in K_0$, then $\gF_{K_0}(x,y)$ belongs to
$S_0$.

Since $\Conc K\iso\two$, there are elements $0_K$ and $1_K$
of $K$ such that $0_K<1_K$ (they are not necessarily the
least and the largest element of $K$, even if the latter
exist), and then the equality $\gF_K(0_K,1_K)=\one$ holds.
Furthermore, put $0_{K_i}=f_i(0_K)$ and $1_{K_i}=f_i(1_K)$, for
all $i<3$. Note, in particular, that the following equalities
hold:
 \[
 \gF_{K_i}(0_{K_i},1_{K_i})=\gF_K(0_K,1_K)=\one.
 \]
Further, put
$0_{L_j}=g_{ij}(0_{K_i})$ and $1_{L_j}=g_{ij}(1_{K_i})$,
for all $j<3$ such that $i\ne j$;
this definition is consistent, for example, the value of
$0_{L_j}=g_{ij}(0_{K_i})$ does not depend of the choice of
$i$, because of the commutativity of the diagram of
Figure~\ref{Fi:DiagLatt}. Finally, put $0_P=h_j(0_{L_j})$
and $1_P=h_j(1_{L_j})$; again, this does not depend of $j$.
\smallskip

For all $i<3$, the equality
$\gF_{K_i}(0_{K_i},1_{K_i})=\one=\ga_i\jj\gb_i$ holds, and
$K_i$ has permutable congruences, thus there exists $x_i\in K_i$
such that $\gF_{K_i}(0_{K_i},x_i)\leq\ga_i$ and
$\gF_{K_i}(x_i,1_{K_i})\leq\gb_i$. By replacing $x_i$ by
$(x_i\jj0_{K_i})\mm1_{K_i}$, one can suppose that
$0_{K_i}\leq x_i\leq1_{K_i}$. If $x_i=0_{K_i}$, then
$\gF_{K_i}(0_{K_i},1_{K_i})\leq\gb_i$, that is,
$\one\leq\gb_i$, a contradiction. Thus, $0_{K_i}<x_i$, thus
$\gF_{K_i}(0_{K_i},x_i)>\zero$. Since
$\gF_{K_i}(0_{K_i},x_i)\leq\ga_i$ and since
$S_i=\set{\zero,\ga_i,\gb_i,\one}$, and by symmetry, we
obtain the following:
\begin{equation}\label{Eq:OblConxi}
\gF_{K_i}(0_{K_i},x_i)=\ga_i\quad\text{and}\quad
\gF_{K_i}(x_i,1_{K_i})=\gb_i.
\end{equation}

By applying to the elements $0_{K_i}$, $x_i$, and $1_{K_i}$
the homomorphisms $g_{ij}$ for $i\ne j$, we obtain
sublattices of the $L_j$ which can be described by
Figure~\ref{Fi:Subl}. We use here the following notation: for
$j<3$, $j'$, and $j''$ are the two elements of
$3\setminus\set{j}$, ordered in such a way that $j'<j''$.

\begin{figure}[htb]
\begin{picture}(100,110)(0,0)
\thicklines
\put(20,0){\circle{6}}
\put(20,30){\circle{6}}
\put(0,50){\circle{6}}
\put(40,50){\circle{6}}
\put(20,70){\circle{6}}
\put(20,100){\circle{6}}

\put(15,0){\makebox(0,0)[r]{$0_{L_j}$}}
\put(15,30){\makebox(0,0)[r]{$u_j\mm v_j$}}
\put(-5,50){\makebox(0,0)[r]{$u_j=g_{j'j}(x_{j'})$}}
\put(45,50){\makebox(0,0)[l]{$v_j=g_{j''j}(x_{j''})$}}
\put(15,70){\makebox(0,0)[r]{$u_j\jj v_j$}}
\put(15,100){\makebox(0,0)[r]{$1_{L_j}$}}

\put(20,3){\line(0,1){24}}
\put(17.88,32.12){\line(-1,1){15.76}}
\put(22.12,32.12){\line(1,1){15.76}}
\put(2.12,52.12){\line(1,1){15.76}}
\put(37.88,52.12){\line(-1,1){15.76}}
\put(20,73){\line(0,1){24}}

\end{picture}
\caption{A sublattice of $L_j$, for $j<3$}
\label{Fi:Subl}
\end{figure}
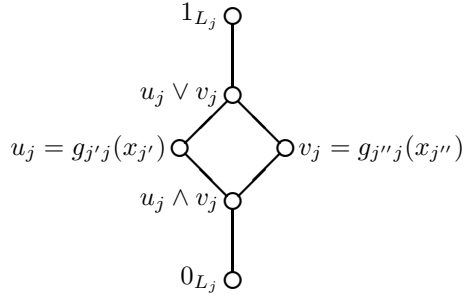

We put $u_j=g_{j'j}(x_{j'})$ and $v_j=g_{j''j}(x_{j''})$,
for all $j<3$.
Note, in particular, that, by \eqref{Eq:OblConxi}, the
following equalities hold:
 \[
 \begin{aligned}
 \gF_{L_j}(0_{L_j},u_j)&=\gF_{K_{j'}}(0_{K_{j'}},x_{j'})
 &&=\ga_{j'},\\
 \gF_{L_j}(u_j,1_{L_j})&=\gF_{K_{j'}}(x_{j'},1_{K_{j'}})
 &&=\gb_{j'}.
 \end{aligned}
 \]
Similarly, the following equalities hold: 
 \[
 \begin{aligned}
 \gF_{L_j}(0_{L_j},v_j)&=\gF_{K_{j''}}(0_{K_{j''}},x_{j''})
 &&=\ga_{j''},\\
 \gF_{L_j}(v_j,1_{L_j})&=\gF_{K_{j''}}(x_{j''},1_{K_{j''}})
 &&=\gb_{j''}.
 \end{aligned}
 \]
Therefore, by applying $\gF_{L_j}$ to the edges of the
graph representing Figure~\ref{Fi:Subl}, we obtain, for all
$j<3$, the following refinement matrix (in $T_j$):
 \begin{equation}\label{Eq:RefMatCon}
 \mbox{\begin{tabular}{|c|c|c|}
 \cline{2-3}
 \multicolumn{1}{l|}{} & \tvi $\ga_{j''}$ & $\gb_{j''}$\\
 \hline
 \tvi $\ga_{j'}$ & $\gF_{L_j}(0_{L_j},u_j\mm v_j)$
 & $\gF_{L_j}(u_j\mm v_j,u_j)$\\
 \hline
 \tvi $\gb_{j'}$ & $\gF_{L_j}(u_j\mm v_j,v_j)$
 & $\gF_{L_j}(u_j\jj v_j,1_{L_j})$\\
 \hline
 \end{tabular}}
 \end{equation}
Since the entries of this matrix lie in $T_j$, and, in
the distributive lattice $T_j$, the following equalities
 \[
 \ga_{j'}\mm\gb_{j'}=\ga_{j''}\mm\gb_{j''}=\zero
 \]
hold, the only possibility is the one given by the
refinement matrices \eqref{Eq:RefMat}.
In particular, we obtain the following relations:
 \begin{equation}\label{Eq:FinalRel}
 \gF_{L_0}(u_0\mm v_0,u_0)=\gx_1,\quad
 \gF_{L_1}(u_1\mm v_1,u_1)=\gh_1,\quad\text{and}\quad
 \gF_{L_2}(u_2\mm v_2,u_2)=\gz_1.
 \end{equation}
Finally, note that the equality
 \begin{equation}\label{Eq:trans}
 \gF_{L_j}(x,y)=\gF_P(h_j(x),h_j(y))
 \end{equation}
holds for all $j<3$ and for all $x$, $y\in L_j$.
Furthermore, the equality
$h_j(u_j)=h_jg_{j'j}(x_{j'})=g_{j'}(x_{j'})$ holds, and,
similarly, $h_j(v_j)=h_jg_{j'j}(x_{j'})=g_{j''}(x_{j''})$.
Thus, by applying \eqref{Eq:trans} to \eqref{Eq:FinalRel},
we obtain the following equalities:
 \begin{align*}
 \gx_1&=\gF_P(g_1(x_1)\mm g_2(x_2),g_1(x_1)),\\
 \gh_1&=\gF_P(g_0(x_0)\mm g_2(x_2),g_0(x_0)),\\
 \gz_1&=\gF_P(g_0(x_0)\mm g_1(x_1),g_0(x_0)).
 \end{align*}
In particular, $\gh_1\leq\gx_1\jj\gz_1$. But this
contradicts Lemma~\ref{L:NonIneq}.
\end{proof}

\begin{remark}
The proof above shows, in fact, a stronger result: namely,
Theorem~\ref{T:NonAmalg} can be generalized to any diagram
of semilattices obtained from $\Dc$ by
replacing $U$ by any \emph{larger} distributive semilattice.
To paraphrase this, the problem for lifting
Figure~\ref{Fi:DiagSemil} lies in the fact that the top
lattice in Figure~\ref{Fi:DiagLatt} would have a too
\emph{small} congruence lattice.
\end{remark}

\section{\bf Duality of complete lattices; case of $\Dc$}
\label{S:Dual}

We shall introduce in this section some material that will
eventually lead to the existence of a lifting of $\Dc$, with
respect to the $\Conc$ functor, by a diagram of finite
lattices, see Theorem~\ref{T:Lift}.

\subsection{Duality of complete lattices}\label{S:DualCpl}

The facts presented in this section are standard, although we do
not know of any reference where they are recorded. All the
proofs are straightforward, so we omit them.

If $A$ and $B$ are complete lattices, a
map $f\colon A\to B$ is a \emph{\cjh}, if the equality
 \[
 f\left(\JJ X\right)=\JJ f[X]
 \]
holds, for every subset $X$ of $A$. Note that this implies, in
particular, that $f(0_A)=0_B$. One defines, similarly,
\emph{\cmh s}. We shall denote by $\Cj$ (resp., $\Cm$) the
category of complete lattices with \cjh s (resp., \cmh s).

\begin{definition}\label{D:Dual}
Let $A$ and $B$ be complete lattices.
Two maps $f\colon A\to B$ and $g\colon B\to A$ are
\emph{dual}, if the equivalence
 \[
 f(a)\leq b\text{ if and only if }a\leq g(b),
 \]
holds, for all $\vv<a,b>\in A\times B$.
\end{definition}

\begin{lemma}\label{L:Dual}
Let $A$ and $B$ be complete lattices.

\begin{enumerate}
\item If $f\colon A\to B$ and $g\colon B\to A$ are dual,
then $f$ is a \cjh\ and $g$ is a \cmh.

\item Let $f\colon A\to B$ be a \cjh. Then there exists a unique
map $g\colon B\to A$ such that $f$ and $g$ are dual.

\item Let $g\colon B\to A$ be a \cmh. Then there exists a unique
map $f\colon A\to B$ such that $f$ and $g$ are dual.

\end{enumerate}

\end{lemma}

In case (ii), for all $b\in B$, $g(b)$ is defined as the largest
$a\in A$ such that $f(a)\leq b$. We let $g$ denote $f^*$.
Similarly, in case (iii), $f(a)$ is defined as the least $b\in B$
such that $a\leq g(b)$, and we let $f$ denote $g^\dagger$.

The basic categorical properties of the duality thus described
may be recorded in the following lemma.

\goodbreak
\begin{lemma}\label{L:CatDual}\hfill
\begin{enumerate}
\item The correspondence $f\mapsto f^*$ defines a contravariant
functor from $\Cj$ to $\Cm$.

\item The correspondence $g\mapsto g^\dagger$ defines a
contravariant functor from $\Cm$ to $\Cj$.

\item If $f$ is a \cjh, then $(f^*)^\dagger=f$.

\item If $g$ is a \cmh, then $(g^\dagger)^*=g$.

\end{enumerate}
\end{lemma}

Of particular importance is the effect of the duality on
\cjh\ of the form $\Con f\colon\Con K\to\Con L$, where
$f\colon K\to L$ is a lattice homomorphism. We denote by
$\Res f\colon\Con L\to\Con K$ the ``restriction'' map, defined by
 \[
 (\Res f)(\gb)=
 \setm{\vv<x,y>\in K\times K}{\vv<f(x),f(y)>\in\gb},
 \]
for all $\gb\in\Con L$. If $f$ is the inclusion mapping from a lattice
$K$ into a lattice $L$, we shall just denote $(\Res f)(\gb)$ by
$\gb\res_K$, and we shall call it the \emph{restriction} of
$\gb$ to~$K$.

\begin{lemma}\label{L:DualExt}
Let $K$ and $L$ be lattices, let $f\colon K\to L$ be a lattice
homomorphism.
Then $\Con f$ and $\Res f$ are dual.
\end{lemma}

For a lattice $L$, we denote by $\bgi_L$ the
largest congruence of $L$.

\begin{lemma}\label{L:DualInt}
Let $K$ and $L$ be lattices, let $f\colon K\to L$ be a lattice
homomorphism, let $\gb\in\Con L$. We define $\gy$ as the
restriction of $\Res f$ to the interval $[\gb,\bgi_L]$ of
$\Con L$, and we put $\ga=\gy(\gb)$.
Furthermore, we denote by $f'\colon K/\ga\to L/\gb$ the
canonical lattice homomorphism, and we put $\gy'=\Res f'$.
Let $\ge\colon[\gb,\bgi_L]\to\Con(L/\gb)$ and
$\gh\colon[\ga,\bgi_K]\to\Con(K/\ga)$ be the canonical
isomorphisms. Then the following diagram commutes:
 \[
 \begin{CD}
 \Con(L/\gb) @>{\gy'}>> \Con(K/\ga)\\
 @A{\ge}AA @AA{\gh}A\\
 [\gb,\bgi_L]Ê@>>{\gy}>  [\ga,\bgi_K]
 \end{CD}
 \]
\end{lemma}

\subsection{The dual of $\Dc$}\label{S:DualDc}

We describe in this section the duals of the semilattice mappings
of the semilattice diagram $\Dc$. Note that we consider join- or
\mhom s between finite lattices, so
these homomorphisms are always complete. The inclusion mappings
in $\Dc$ are \jz-homomorphisms, but they are not always
\mhom s. Since the duals of these inclusion mappings
are \mhom s, we only need to specify their values at
the meet-irreducible elements (that is, the coatoms) of the
Boolean semilattices belonging to $\Dc$.

The coatoms of $U$ are the subsets $\bar k=5\setminus\set{k}$,
for all $k<5$. Now, the Boolean subsemilattice $T_0$ of $U$ is
generated by the atoms $\gx_0=\set{0,4}$, $\gx_1=\set{3}$,
$\gx_2=\set{2}$, $\gx_3=\set{1,4}$. Each coatom of $T_0$ is also
a coatom of $U$. So the coatoms of $T_0$ are
 \[
 \bar{\gx}_0=\gx_1\jj\gx_2\jj\gx_3=\set{1,2,3,4}=\bar 0,\quad
 \bar{\gx}_1=\bar 3,\quad\bar{\gx}_2=\bar 2,\quad
 \bar{\gx}_3=\bar 1.
 \]
We obtain, similarly, the coatoms of $T_1$, respectively $T_2$:
 \begin{align*}
 \bar{\gh}_0&=\bar 0,&\bar{\gh}_1&=\bar 1,&
 \bar{\gh}_2&=\bar 2,&\bar{\gh}_3&=\bar 3,\\
 \bar{\gz}_0&=\bar 0,&\bar{\gz}_1&=\bar 1,&
 \bar{\gz}_2&=\bar 3,&\bar{\gz}_3&=\bar 2.
 \end{align*}
The atoms of the Boolean semilattices $S_i$, for $i<3$, are also
their coatoms, so we still denote them by $\ga_i$, $\gb_i$,
for $i<3$. The bottom lattice, $\two$, has a unique coatom,
namely, $\zero$.

Next, we describe the mappings $\gy_{j,i}\colon T_j\to S_i$ that
are the duals of the inclusion maps $S_i\hookrightarrow T_j$, for
$i\ne j$. Easy computations yield the following:
\begin{equation}\label{Eq:psimess}
 \begin{aligned}
 \psi_{0,2}(\bar 1)&=\psi_{0,2}(\bar 3)=\alpha_2,&\qquad
 \psi_{0,2}(\bar 0)&=\psi_{0,2}(\bar 2)=\beta_2,\\
 \psi_{0,1}(\bar 1)&=\psi_{0,1}(\bar 2)=\alpha_1,&\qquad
 \psi_{0,1}(\bar 0)&=\psi_{0,1}(\bar 3)=\beta_1,\\
 \psi_{1,2}(\bar 1)&=\psi_{1,2}(\bar 3)=\alpha_2,&\qquad
 \psi_{1,2}(\bar 0)&=\psi_{1,2}(\bar 2)=\beta_2,\\
 \psi_{1,0}(\bar 2)&=\psi_{1,0}(\bar 3)=\alpha_0,&\qquad
 \psi_{1,0}(\bar 0)&=\psi_{1,0}(\bar 1)=\beta_0,\\
 \psi_{2,1}(\bar 1)&=\psi_{2,1}(\bar 2)=\alpha_1,&\qquad
 \psi_{2,1}(\bar 0)&=\psi_{2,1}(\bar 3)=\beta_1,\\
 \psi_{2,0}(\bar 2)&=\psi_{2,0}(\bar 3)=\alpha_0,&\qquad
 \psi_{2,0}(\bar 0)&=\psi_{2,0}(\bar 1)=\beta_0.
 \end{aligned}
\end{equation}
The computations of the values are always based on the fact that each
of the (co)atoms $\ga_i$, $\gb_i$, for $i<3$, omits exactly two
elements of the set $\set{0,1,2,3}$ and hence it lies under exactly
two coatoms of each semilattice $T_j$, for $j<3$.

Next we describe the mappings $\gy_j\colon U\to T_j$, for
$j<3$, the duals of the inclusion mappings
$T_j\hookrightarrow U$. It is easy to describe the values of
$\gy_j$ at the coatoms $\bar 0$, $\bar 1$, $\bar 2$,
$\bar 3\in U$, since each of them is also a coatom of $T_j$,
thus
 \[
 \gy_j(\bar 0)=\bar 0,\quad\gy_j(\bar 1)=\bar 1,\quad
 \gy_j(\bar 2)=\bar 2,\quad\gy_j(\bar 3)=\bar 3,
 \]
for all $j<3$. The values $\gy_j(\bar 4)$ are crucial.
So $\gy_0(\bar 4)$ is the largest element of $T_0$ not containing
$4$ as an element, hence
 \[
 \gy_0(\bar 4)=\gx_1\jj\gx_2=\set{2,3}=\bar 0\mm\bar 1,
 \]
where the meet is computed in $T_0$, of course. Similarly, we obtain
that
 \[
 \begin{aligned}
 \gy_1(\bar 4)&=\set{2}&&=\bar 0\mm\bar 1\mm\bar 3,\\
 \gy_2(\bar 4)&=\set{1,3}&&=\bar 0\mm\bar 2.
 \end{aligned}
 \]

Thus we also get $\gy_{j,i}\gy_j(\bar 4)=\ga_i\mm\gb_i=\zero$,
for all $i\ne j$. Moreover, for any other coatom $\bar k$, $k<4$,
we get that $\gy_{j,i}\gy_j(\bar k)$ is the unique of the
(co)atoms $\ga_i$, $\gb_i$ not containing $k$ as an element. It
is computed by the same formulas as those in \eqref{Eq:psimess}.

The dual of any inclusion mapping $\two\hookrightarrow X$,
where $X$ is any one of the Boolean lattices $U$, $T_j$, $S_i$,
maps every coatom of
$X$ to the unique coatom of $\two$, namely to~$\zero$.

\section{A lifting of $\Dc$}\label{S:liftDc}

We shall construct in this section a lifting of $\Dc$, based on
the computations of Section~\ref{S:DualDc}. The computations
in this section are relatively tedious; however, in our opinion,
they carry the hope of being generalizable to further situations,
like being able to lift arbitrary $\two^n$-diagrams of
finite distributive \jz-semilattices. This purpose in mind, we
found it useful to give the computations in some detail. 

It is also important to note that the lattices that constitute
our lifting do \emph{not} have permutable congruences---they
cannot, by Theorem~\ref{T:NonAmalg}.\smallskip

For a variety $\VV$ of lattices and a set $X$, we denote by
$\FV(X)$ the free lattice in $\VV$ generated by $X$. If
$Y\ci X$, then $\FV(Y)$ is a sublattice of $\FV(X)$. In fact,
$\FV(Y)$ is a \emph{retract} of $\FV(X)$. Indeed, let $f$ be any
map from $X$ to $\FV(Y)$ such that $f\res_Y=\id_Y$. We still
denote by $f$ the unique lattice homomorphism from $\FV(X)$ to
$\FV(Y)$ that extends $f$. If $j$ denotes the inclusion map
from $\FV(Y)$ into $\FV(X)$, then $f\circ j=\id_{\FV(Y)}$, which
proves our claim. In particular, the equality
 \[
 (\Con f)\circ(\Con j)=\id_{\Con\FV(Y)}
 \]
holds, thus $\Con j$ is one-to-one, which means that $j$ has the
congruence extension property.

We shall fix in this section a seven-element set,
 \[
 X=\set{a,b,c,d,e,u,v},
 \]
and we define subsets $Y$, $Y_i$, $X_i$ (for $i<3$) of $X$ by
 \begin{align*}
 X_0&=\set{a,b,c,e}, & X_1&=\set{a,b,c,d}, & X_2&=\set{a,b,d,e}\\
 Y_0&=\set{a,b,d},   & Y_1&=\set{a,b,e},   & Y_2&=\set{a,b,c},\\
 &&Y&=\set{a,b}.&&
 \end{align*}

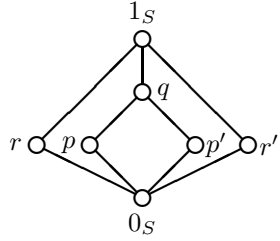
\begin{figure}[htb]
\begin{picture}(100,70)(0,0)
\thicklines

\put(40,0){\circle{6}}
\put(0,20){\circle{6}}
\put(20,20){\circle{6}}
\put(60,20){\circle{6}}
\put(80,20){\circle{6}}
\put(40,40){\circle{6}}
\put(40,60){\circle{6}}

\put(68,20){\makebox(0,0){$p'$}}
\put(48,40){\makebox(0,0){$q$}}
\put(88,20){\makebox(0,0){$r'$}}
\put(12,20){\makebox(0,0){$p$}}
\put(-8,20){\makebox(0,0){$r$}}
\put(40,-10){\makebox(0,0){$0_S$}}
\put(40,70){\makebox(0,0){$1_S$}}

\put(37.32,1.34){\line(-2,1){34.63}}
\put(42.68,1.34){\line(2,1){34.63}}
\put(37.88,2.12){\line(-1,1){15.76}}
\put(42.12,2.12){\line(1,1){15.76}}
\put(22.12,22.12){\line(1,1){15.76}}
\put(57.88,22.12){\line(-1,1){15.76}}
\put(40,43){\line(0,1){14}}
\put(2.12,22.12){\line(1,1){35.76}}
\put(77.88,22.12){\line(-1,1){35.76}}

\end{picture}
\caption{The lattice $S$}
\label{Fi:SimpLat}
\end{figure}

Let $S$ be the lattice diagrammed on Figure~\ref{Fi:SimpLat}.

Note, in particular, that $S$ is a finite, simple, non-modular
lattice. In fact, the argument could be carried out for
\emph{any} simple, bounded, non-modular lattice instead of $S$. We
denote by $\VV$ the variety generated by $S$.

We define lattices $A_k$, for $k<5$, as follows:
 \[
 A_0=A_1=A_2=A_3=\two,\quad A_4=S.
 \]
Note that all the $A_k$-s are simple lattices of $\VV$. Next, we
define maps $f_k\colon X\to A_k$, for $k<5$, as follows. For
$k<4$, the $f_k$-s are uniquely determined by
 \begin{equation}\label{Eq:fkXk<4}
 \begin{aligned}
 f_0(a)&=f_0(e)&&<f_0(b)&&=f_0(c)&&=f_0(d)&&=f_0(u)&&=f_0(v),\\
 f_1(a)&=f_1(c)&&<f_1(b)&&=f_1(d)&&=f_1(e)&&=f_1(u)&&=f_1(v),\\
 f_2(a)&=f_2(d)&&<f_2(b)&&=f_2(c)&&=f_2(e)&&=f_2(u)&&=f_2(v),\\
 f_3(a)&=f_3(c)&&=f_3(d)&&=f_3(e)&&<f_3(b)&&=f_3(u)&&=f_3(v),\\
 \end{aligned}
 \end{equation}
while $f_4$ is determined by
 \begin{equation}\label{Eq:f4X}
 \begin{aligned}
 f_4(a)&=0_S,&\q f_4(b)&=1_S,&\q f_4(c)&=p, &&\\
 f_4(d)&=q,  &\q  f_4(e)&=r,  &\q f_4(u)&=p',&\q f_4(v)&=r'.
 \end{aligned}
 \end{equation}
Note that $f_k[X]$ generates $A_k$, for all $k<5$. Since $A_k$
belongs to $\VV$, $f_k$ induces a unique surjective lattice
homomorphism from $\FV(X)$ onto $A_k$, that we still denote by
$f_k$.

We put $\gr_k=\ker f_k$, the kernel of $f_k$, for all $k<5$.
Since $A_k$ is a simple lattice, $\gr_k$ is a coatom of
$\Con\FV(X)$. Furthermore, for $i\ne j$, $f_i$ and $f_j$ define
distinct partitions of $X$, thus $\gr_i\ne\gr_j$. We define a
congruence $\gQ$ of $\FV(X)$ and a lattice $P$, by
 \begin{equation}\label{Eq:gQintgrk}
 \gQ=\MMm{\gr_k}{k<5},\quad P=\FV(X)/\gQ.
 \end{equation}
The congruence lattice of $P$ can be easily computed by using the
following folklore lemma:

\begin{lemma}\label{L:FolDis}
Let $D$ be a distributive lattice with unit, let $n\in\go$, let
$a_i$, for $i<n$, be mutually distinct coatoms of $D$.
We put $a=\MMm{a_i}{i<n}$. Then the
interval $[a,1]$ of $D$ is isomorphic to $\two^n$.
\end{lemma}

By using Lemma~\ref{L:FolDis} for $D=\Con\FV(X)$, $n=5$, and
$a_k=\gr_k$, for $k<5$, we obtain that the upper interval
$[\gQ,\bgi_{\FV(X)}]$ of $\Con\FV(X)$ is isomorphic to $\two^5$.
Thus, $\Con P\iso\two^5$. The coatoms of $\Con P$ are the
congruences $\gr_k/\gQ$, for $k<5$.
Therefore, we have obtained:

\begin{lemma}\label{L:ConP}
$\Con P\iso\two^5$, and the elements $\gr_{k}/\gQ$, for $k<5$,
are the distinct coatoms of $\Con P$.
\end{lemma}

We now construct a commutative cube of lattices, as on
Figure~\ref{Fi:DiagLatt}. For $i<3$, we denote by $\gQ_i$ (resp.,
$\gQ'_i$) the restriction of $\gQ$ to $\FV(X_i)$ (resp., to
$\FV(Y_i)$).
Furthermore, we denote by $\gF$ the restriction of
$\gQ$ to $\FV(Y)$.
We define lattices $K$, $K_i$, and $L_i$, for $i<3$, by
 \[
 K=\FV(Y)/\gF,\quad K_i=\FV(Y_i)/\gQ'_i,\quad
 L_i=\FV(X_i)/\gQ_i.
 \]
Since the relations
 \[
 Y\ci Y_i\ci X_j\ci X
 \]
hold for all $i\ne j$ in $\set{0,1,2}$, there are induced
lattice embeddings $f_i\colon K\hookrightarrow K_i$,
$g_{i,j}\colon K_i\hookrightarrow L_j$, and
$h_j\colon L_j\hookrightarrow P$, for all
$i\ne j$ in $\set{0,1,2}$. It is obvious that these maps form a
commutative diagram as on Figure~\ref{Fi:DiagLatt}. We denote by
\LL\ this diagram. The rest of this section is devoted to the
proof of the following result:

\begin{theorem}\label{T:Lift}
The image of the diagram \LL\ under the $\Con$ functor is
isomorphic to $\Dc$.
\end{theorem}

In order to prove Theorem~\ref{T:Lift}, if is sufficient, by
Section~\ref{S:DualCpl}, to prove that the image under the $\Res$
functor of the diagram \LL\ is isomorphic to the dual diagram of
$\Dc$, described by the maps $\gy_j\colon U\to T_j$,
$\gy_{j,i}\colon T_j\to S_i$, and $\gf_i\colon S_i\to\two$, for
$i\ne j$ in $\set{0,1,2}$. Such an isomorphism of diagrams would
consist of a family of eight isomorphisms, respectively from
$\Con P$ onto
$U$, from $\Con L_i$ onto $T_i$, from $\Con K_i$ onto $S_i$ (for
$i<3$), and from $\Con K$ onto $\two$, satisfying a certain set of
twelve commutation relations, each of them to be verified on the
corresponding set of at most five meet-irreducible elements.
Instead of writing down those cumbersome
relations, it is convenient to observe that the dual of
$\Dc$ may be described by the following data:

\begin{enumerate}
\item $U$ is Boolean, and it has the coatoms $\bar k$, for $k<5$.

\item $\gy_i(\bar k)=\bar k$, a coatom of $T_j$, for all $i<3$
and all $k<4$.

\item $\gy_0(\bar 4)=\bar 0\mm\bar 1$,
$\gy_1(\bar 4)=\bar 0\mm\bar 1\mm\bar 3$,
$\gy_2(\bar 4)=\bar 0\mm\bar 2$.

\item $T_0$, $T_1$, and $T_2$ are Boolean, and they have
the coatoms $\bar k$, for $k<4$.

\item The equations \eqref{Eq:psimess}.

\item $S_i$ is Boolean, and $\ga_i$, $\gb_i$ are the coatoms of
$S_i$, for all $i<3$.

\item For $i<3$, $\gf_i$ is the map that sends $\one$ to
$\one$, and all the other elements to $\zero$.
\end{enumerate}

Let us first highlight which elements of $\Con P$, $\Con L_j$, and
$\Con K_i$, will play the role of the $\bar k$, for $k<5$, and
the $\ga_i$, $\gb_i$, for $i<3$.

We have already seen in Lemma~\ref{L:ConP} that the coatoms of
$\Con P$ are the $\gr_k/\gQ$, for $k<5$. Let $\gr_k/\gQ$
correspond to $\bar k$, for $k<5$.

Now the candidates for the coatoms of $L_j$, for $j<3$. Let
$f_{j,k}\colon\FV(X_j)\to A_k$ be the restriction of $f_k$ to
$\FV(X_j)$, for $k<5$. We put $\gr_{j,k}=\ker f_{j,k}$. So
$\gr_{j,k}$ is the restriction of $\gr_k$ to $\FV(X_j)$.

\begin{lemma}\label{L:ConLj}
Let $j<3$. Then $\Con L_j\iso\two^4$, and the elements
$\gr_{j,k}/\gQ_j$, for $k<4$, are the distinct coatoms of
$\Con L_j$. Furthermore, the following equalities hold:
 \begin{equation}\label{Eq:grj4}
 \begin{aligned}
 \gr_{0,4}&=\gr_{0,0}\mm\gr_{0,1},\\
 \gr_{1,4}&=\gr_{1,0}\mm\gr_{1,1}\mm\gr_{1,3},\\
 \gr_{2,4}&=\gr_{2,0}\mm\gr_{2,2}.
 \end{aligned}
 \end{equation}
\end{lemma}

\begin{proof}
Let us first verify this for $L_0$. We recall that
$X_0=\set{a,b,c,e}$. We compute the values of the maps $f_{0,k}$
on the elements of $X_0$, by just looking at \eqref{Eq:fkXk<4}
and \eqref{Eq:f4X} and removing all the elements of
$X\setminus X_0$:
 \begin{equation}\label{Eq:f0kXk<4}
 \begin{aligned}
 f_{0,0}(a)&=f_{0,0}(e)&&<f_{0,0}(b)&&=f_{0,0}(c),\\
 f_{0,1}(a)&=f_{0,1}(c)&&<f_{0,1}(b)&&=f_{0,1}(e),\\
 f_{0,2}(a)&<f_{0,2}(b)&&=f_{0,2}(c)&&=f_{0,2}(e),\\
 f_{0,3}(a)&=f_{0,3}(c)&&=f_{0,3}(e)&&<f_{0,3}(b),\\
 \end{aligned}
 \end{equation}
while $f_{0,4}$ is determined by
 \begin{equation}\label{Eq:f04X}
 f_{0,4}(a)=0_S,\quad f_{0,4}(b)=1_S,
 \quad f_{0,4}(c)=p,\quad f_{0,4}(e)=r.
 \end{equation}
Since $f_{0,k}[X_0]$ generates the range of $f_{0,k}$, the range
of $f_{0,k}$ equals $\two$, if $k<4$, and
$A_{0,4}=\set{0_S,1_S,p,r}$, a two-atom Boolean lattice, if $k=4$.
In particular, for $k<4$,
$\gr_{0,k}$ is a coatom of $\Con\FV(X_0)$. These congruences are
mutually distinct, because, by \eqref{Eq:f0kXk<4}, they induce
different partitions of $X_0$. Furthermore, we consider the
natural projections $\gx_0$, $\gx_1$ from $A_{0,4}$ to $\two$,
defined by $\gx_0(r)=\gx_1(p)=0$ and $\gx_0(p)=\gx_1(r)=1$. By
checking on the elements of $X_0$, we obtain easily that
$f_{0,0}=\gx_0\circ f_{0,4}$ and $f_{0,1}=\gx_1\circ f_{0,4}$.
Since the map $x\mapsto\vv<\gx_0(x),\gx_1(x)>$ from $A_{0,4}$ to
$\two^2$ is a lattice embedding, the
equation $\gr_{0,4}=\gr_{0,0}\mm\gr_{0,1}$ follows.

We do the same for $L_1$. We recall that
$X_1=\set{a,b,c,d}$. As in the previous paragraph, we compute the
values of the maps $f_{1,k}$ on the elements of $X_1$:
 \begin{equation}\label{Eq:f1kXk<4}
 \begin{aligned}
 f_{1,0}(a)&<f_{1,0}(b)&&=f_{1,0}(c)&&=f_{1,0}(d),\\
 f_{1,1}(a)&=f_{1,1}(c)&&<f_{1,1}(b)&&=f_{1,1}(d),\\
 f_{1,2}(a)&=f_{1,2}(d)&&<f_{1,2}(b)&&=f_{1,2}(c),\\
 f_{1,3}(a)&=f_{1,3}(c)&&=f_{1,3}(d)&&<f_{1,3}(b),\\
 \end{aligned}
 \end{equation}
while $f_{1,4}$ is determined by
 \begin{equation}\label{Eq:f14X}
 f_{1,4}(a)=0_S,\quad f_{1,4}(b)=1_S,
 \quad f_{1,4}(c)=p,\quad f_{1,4}(d)=q.
 \end{equation}
So the range of $f_{1,k}$ equals $\two$, if $k<4$, and
$A_{1,4}=\set{0_S,1_S,p,q}$, a four-element chain,
if $k=4$. In particular, for $k<4$,
$\gr_{1,k}$ is a coatom of $\Con\FV(X_1)$. These congruences are
mutually distinct. Furthermore, we consider the
natural projections $\gh_0$, $\gh_1$, and $\gh_2$ from $A_{1,4}$
to $\two$, defined by $\gh_0(p)=1$,
$\gh_1(p)=0$, $\gh_1(q)=1$, and
$\gh_2(p)=\gh_2(q)=0$, $\gh_2(1)=1$.
By checking on the elements of $X_1$, we obtain easily that
$f_{1,0}=\gh_0\circ f_{1,4}$, $f_{1,1}=\gh_1\circ f_{1,4}$,
and $f_{1,3}=\gh_2\circ f_{1,4}$.
Since the map $x\mapsto\vv<\gh_0(x),\gh_1(x),\gh_2(x)>$ from
$A_{1,4}$ to $\two^3$ is a lattice embedding, the
equation $\gr_{1,4}=\gr_{1,0}\mm\gr_{1,1}\mm\gr_{1,3}$ follows.

We do it finally for $L_2$. We recall that
$X_2=\set{a,b,d,e}$. We compute the
values of the maps $f_{2,k}$ on the elements of $X_2$:
 \begin{equation}\label{Eq:f2kXk<4}
 \begin{aligned}
 f_{2,0}(a)&=f_{2,0}(e)&&<f_{2,0}(b)&&=f_{2,0}(d),\\
 f_{2,1}(a)&<f_{2,1}(b)&&=f_{2,1}(d)&&=f_{2,1}(e),\\
 f_{2,2}(a)&=f_{2,2}(d)&&<f_{2,2}(b)&&=f_{2,2}(e),\\
 f_{2,3}(a)&=f_{2,3}(d)&&=f_{2,3}(e)&&<f_{2,3}(b),\\
 \end{aligned}
 \end{equation}
while $f_{1,4}$ is determined by
 \begin{equation}\label{Eq:f24X}
 f_{2,4}(a)=0_S,\quad f_{2,4}(b)=1_S,
 \quad f_{2,4}(d)=q,\quad f_{2,4}(e)=r.
 \end{equation}
So the range of $f_{2,k}$ equals $\two$, if $k<4$, and
$A_{2,4}=\set{0_S,1_S,q,r}$, a two-atom Boolean lattice,
if $k=4$. In particular, for $k<4$,
$\gr_{2,k}$ is a coatom of $\Con\FV(X_2)$. These congruences are
mutually distinct. Furthermore, we consider the
natural projections $\gz_0$ and $\gz_1$ from $A_{2,4}$
to $\two$, defined by $\gz_0(r)=\gz_1(q)=0$,
and $\gz_0(q)=\gz_1(r)=1$.
By checking on the elements of $X_2$, we obtain easily that
$f_{2,0}=\gz_0\circ f_{2,4}$ and $f_{2,2}=\gz_1\circ f_{2,4}$.
Since the map $x\mapsto\vv<\gz_0(x),\gz_1(x)>$ from
$A_{2,4}$ to $\two^2$ is a lattice embedding, the
equation $\gr_{2,4}=\gr_{2,0}\mm\gr_{2,2}$ follows.

In particular, it follows from \eqref{Eq:gQintgrk},
\eqref{Eq:grj4} that the equation
 \[
 \gQ_j=\MMm{\gr_{j,k}}{k<4}
 \]
holds for all $j<3$. Thus the $\gr_{j,k}/\gQ_j$, for $k<4$, are
exactly the coatoms of $L_j$.
\end{proof}

At this point, we have verified (i)--(iv) of the data
that describe the dual of $\Dc$:

\begin{itemize}
\item[(i)] $\Con P$ is Boolean, and it has the coatoms
$\tgr_k=\gr_k/\gQ$, for $k<5$.

\item[(ii)] We put $\tgr_{j,k}=(\Res h_j)(\tgr_k)$,
for all $j<3$ and all $k<4$. By Lemma~\ref{L:DualInt},
$\tgr_{j,k}=\gr_{j,k}/\gQ_j$. If $k<4$, then $\tgr_{j,k}$ is a
coatom of $\Con L_j$.

\item[(iii)] The following equations hold:
 \begin{equation}\label{Eq:tgrj4}
 \begin{aligned}
 (\Res h_0)(\tgr_4)=\tgr_{0,4}&=\tgr_{0,0}\mm\tgr_{0,1},\\
 (\Res h_1)(\tgr_4)=\tgr_{1,4}
 &=\tgr_{1,0}\mm\tgr_{1,1}\mm\tgr_{1,3},\\
 (\Res h_2)(\tgr_4)=\tgr_{2,4}&=\tgr_{2,0}\mm\tgr_{2,2}.
 \end{aligned}
 \end{equation}

\item[(iv)] By Lemma~\ref{L:ConLj}, $\Con L_j$ is Boolean, for
$j<3$, and its coatoms are the $\tgr_{j,k}$-s, for $k<4$.

\end{itemize}

We proceed through the verification of (v)--(vii). The analogues
of $\ga_i$, $\gb_i$ have not been defined yet. We do this now.

For $i<3$ and $k<5$, we denote by $f'_{i,k}$ the restriction of
$f_k$ to $\FV(Y_i)$, and by $\gr'_{i,k}$ the kernel of
$f'_{i,k}$. If $j\ne i$ in $\set{0,1,2}$, then $\gr'_{i,k}$ is
the restriction of $\gr_{j,k}$ to $\FV(X_i)$. In particular, by
the equations \eqref{Eq:grj4} in Lemma~\ref{L:ConLj},
$\gr'_{i,4}$ is the meet of elements of the form $\gr'_{i,k}$,
for $k<4$. Thus, in order to determine the meet-irreducible
elements of $\Con K_i$, it is sufficient to compute $\gr'_{i,k}$,
for $k<4$. We follow a similar, though slightly simpler, pattern
as in the proof of Lemma~\ref{L:ConLj}.

We first compute $f'_{0,k}$, for $k<4$, at the elements of
$Y_0=\set{a,b,d}$, by using \eqref{Eq:fkXk<4}.
We obtain the following:
 \begin{equation}\label{Eq:f'0kXk<4}
 \begin{aligned}
 f'_{0,0}(a)&<f'_{0,0}(b)&&=f'_{0,0}(d),\\
 f'_{0,1}(a)&<f'_{0,1}(b)&&=f'_{0,1}(d),\\
 f'_{0,2}(a)&=f'_{0,2}(d)&&<f'_{0,2}(b),\\
 f'_{0,3}(a)&=f'_{0,3}(d)&&<f'_{0,3}(b).\\
 \end{aligned}
 \end{equation}
In particular, $f'_{0,0}=f'_{0,1}$ and $f'_{0,2}=f'_{0,3}$.
So we put $\bga_0=\gr'_{0,2}=\gr'_{0,3}$, and
$\bgb_0=\gr'_{0,0}=\gr'_{0,1}$. We note that $\bga_0\ne\bgb_0$,
and that the range of $f'_{0,k}$ is isomorphic to $\two$, for
$k<4$. So, $\bga_0$ and $\bgb_0$ are distinct coatoms of
$\Con\FV(Y_0)$. By \eqref{Eq:gQintgrk}, they meet to $\gQ'_0$.
Hence, $\Con K_0\iso\two^2$, and the coatoms of $\Con K_0$ are
$\tga_0=\bga_0/\gQ'_0$ and $\tgb_0=\bgb_0/\gQ'_0$.

Similarly, we compute $f'_{1,k}$, for $k<4$, at the elements of
$Y_1=\set{a,b,e}$. We obtain the following:
 \begin{equation}\label{Eq:f'1kXk<4}
 \begin{aligned}
 f'_{1,0}(a)&=f'_{1,0}(e)&&<f'_{1,0}(b),\\
 f'_{1,1}(a)&<f'_{1,1}(b)&&=f'_{1,1}(e),\\
 f'_{1,2}(a)&<f'_{1,2}(b)&&=f'_{1,2}(e),\\
 f'_{1,3}(a)&=f'_{1,3}(e)&&<f'_{1,3}(b).\\
 \end{aligned}
 \end{equation}
In particular, $f'_{1,0}=f'_{1,3}$ and $f'_{1,1}=f'_{1,2}$.
So we put $\bga_1=\gr'_{1,1}=\gr'_{1,2}$, and
$\bgb_1=\gr'_{1,0}=\gr'_{1,3}$. We note that $\bga_1\ne\bgb_1$,
and that the range of $f'_{1,k}$ is isomorphic to~$\two$, for
$k<4$. So, $\bga_1$ and $\bgb_1$ are distinct coatoms of
$\Con\FV(Y_1)$. They meet to $\gQ'_1$.
Hence, $\Con K_1\iso\two^2$, and the coatoms of $\Con K_1$ are
$\tga_1=\bga_1/\gQ'_1$ and $\tgb_1=\bgb_1/\gQ'_1$.

Finally, we compute $f'_{2,k}$, for $k<4$, at the elements of
$Y_2=\set{a,b,c}$. We obtain the following:
 \begin{equation}\label{Eq:f'2kXk<4}
 \begin{aligned}
 f'_{2,0}(a)&<f'_{2,0}(b)&&=f'_{2,0}(c),\\
 f'_{2,1}(a)&=f'_{2,1}(c)&&<f'_{2,1}(b),\\
 f'_{2,2}(a)&<f'_{2,2}(b)&&=f'_{2,2}(c),\\
 f'_{2,3}(a)&=f'_{2,3}(c)&&<f'_{2,3}(b).\\
 \end{aligned}
 \end{equation}
In particular, $f'_{2,0}=f'_{2,2}$ and $f'_{2,1}=f'_{2,3}$.
So we put $\bga_2=\gr'_{2,1}=\gr'_{2,3}$, and
$\bgb_2=\gr'_{2,0}=\gr'_{2,2}$. We note that $\bga_2\ne\bgb_2$,
and that the range of $f'_{2,k}$ is isomorphic to~$\two$, for
$k<4$. So, $\bga_2$ and $\bgb_2$ are distinct coatoms of
$\Con\FV(Y_2)$. They meet to $\gQ'_2$.
Hence, $\Con K_2\iso\two^2$, and the coatoms of $\Con K_2$ are
$\tga_2=\bga_2/\gQ'_2$ and $\tgb_2=\bgb_2/\gQ'_2$.

This takes care of (vi) of the data describing the dual of $\Dc$:
by \eqref{Eq:f'0kXk<4}--\eqref{Eq:f'2kXk<4}, $\Con K_i$ is
Boolean and has distinct coatoms $\tga_i$, $\tgb_i$, for all
$i<3$.

Now we verify (v). We just do the typical case from $\Con L_0$ to
$\Con K_2$, the other five proceeding in a similar fashion. The
computation is, actually, easy:
 \[
 \begin{aligned}
 (\Res g_{2,0})(\tgr_{0,0})&=\gr'_{2,0}/\gQ'_2&&=\tgb_2,\\
 (\Res g_{2,0})(\tgr_{0,1})&=\gr'_{2,1}/\gQ'_2&&=\tga_2,\\
 (\Res g_{2,0})(\tgr_{0,2})&=\gr'_{2,2}/\gQ'_2&&=\tgb_2,\\
 (\Res g_{2,0})(\tgr_{0,3})&=\gr'_{2,3}/\gQ'_2&&=\tga_2.
 \end{aligned}
 \]
Hence, $\Res g_{2,0}$ acts on $\tgr_{0,k}$ as the
map $\gy_{0,2}$ acts on $\bar k$, for $k<4$. Similarly, we can
prove that for $i\ne j$ in $\set{0,1,2}$,
$\Res g_{j,i}$ acts on $\tgr_{i,k}$, for $k<4$, as the
map $\gy_{i,j}$ acts on $\bar k$, for $k<4$.

The verification of (vi) is easy. Since $f_k(a)<f_k(b)$ for
all $k<5$, the restriction mapping from every interval
$[\gQ'_i,\bgi_{\FV(Y_i)}]$, for $i<3$, maps every coatom to
the only coatom of $\two$, namely, $\zero$. Thus $\Res f_i$
lifts the dual $\gf_i$ of the inclusion mapping
$\two\hookrightarrow S_i$.

This completes the proof of Theorem~\ref{T:Lift}. Note that
since $S$ is a finite lattice, the variety $\VV$ is locally
finite, so all lattices $\FV(X)$, $\FV(X_i)$, $\FV(Y_i)$, and
$\FV(Y)$, for $i<3$, are finite. Hence, \emph{a fortiori}, all
lattices $P$, $L_i$, $K_i$, and $K$, for $i<3$, are finite. This
proves that \emph{the diagram $\Dc$ has a lifting by finite
lattices and lattice homomorphisms}. It is, in fact, easy to
prove that $K_i$, for $i<3$, is a three-element chain. In
particular, $K_i$ has almost permutable congruences.

\section{A cube of finite Boolean semilattices without a lifting
by lattices with almost permutable congruences}\label{S:cubeac}

We construct in this section an extension, $\Dac$, of the
semilattice cube $\Dc$ described in Section~\ref{S:SemCube},
that cannot be lifted by lattices with \emph{almost} permutable
congruences. This gives a combinatorial analogue of
Corollary~\ref{C:FreeNonType1}.

The finite semilattices in the cube will again be
subsemilattices of a Boolean lattice, this time on $8$ elements.
So $U$ is, this time, the semilattice of all subsets of the set
$8=\set{0,1,2,3,4,5,6,7}$. We define elements
$\gx_i$, $\gh_i$, $\gz_i$, for $i\in\set{0,1,2,3}$, as follows:
 \begin{align*}
 \gx_0&=\set{0,4,7},&\gx_1&=\set{3,5,6},
 &\gx_2&=\set{2,5,6},&\gx_3&=\set{1,4,7};\\
 \gh_0&=\set{0,4,5,7},&\gh_1&=\set{1,4,6,7},
 &\gh_2&=\set{2,5,6,7},&\gh_3&=\set{3,4,5,6};\\
 \gz_0&=\set{0,4,6},&\gz_1&=\set{1,5,7},
 &\gz_2&=\set{3,5,7},&\gz_3&=\set{2,4,6}.
 \end{align*}
We denote by $T_0$ the \jz-subsemilattice generated by
$\setm{\gx_j}{j<4}$. Because of the elements
$0$, $1$, $2$, $3$, $T_0$
is isomorphic to the Boolean semilattice of all subsets of a
four-element set. Similarly, the \jz-subsemilattice
$T_1$ of $U$ generated by $\setm{\gh_j}{j<4}$ is also
isomorphic to the semilattice of all subsets of a four-element
set. Similarly, the \jz-subsemilattice $T_2$ of $U$ generated
by $\setm{\gz_j}{j<4}$ is isomorphic to both $T_0$ and $T_1$.

Further, we denote by $\ga_i$, $\gb_i$, $i<3$, the following subsets of
$\set{0,1,2,3,4,5,6,7}$:
 \begin{align*}
 \ga_0&=\set{0,1,4,5,6,7},&\gb_0&=\set{2,3,4,5,6,7};\\
 \ga_1&=\set{0,3,4,5,6,7},&\gb_1&=\set{1,2,4,5,6,7};\\
 \ga_2&=\set{0,2,4,5,6,7},&\gb_2&=\set{1,3,4,5,6,7}.
 \end{align*}
If $S_i$, for $i<3$, denotes the \jz-subsemilattice of $U$
generated by $\ga_i$, $\gb_i$, then each semilattice $S_i$ is
isomorphic to $\two^2$. Moreover, $S_i\subseteq T_j$ if $i\ne j$.
The bottom semilattice of the cube is $\two=\set{\es,8}$.

We denote by $\Dac$ this new diagram of finite Boolean
\jz-semilattices. It has the same shape as $\Dc$, but has new
values for the $S_i$, $T_j$, and $U$.

\begin{theorem}\label{T:NonAmalg2}
There exists no lifting, with respect to the $\Conc$ functor,
in the category of lattices, of the diagram $\Dac$, such that
the lattices corresponding to $S_i$, for $i<3$, have almost
permutable congruences.
\end{theorem}

\begin{proof}[Outline of proof]
We merely outline the proof here, by indicating the modifications that
have to be performed on the proof of Theorem~\ref{T:NonAmalg}.
Lemma~\ref{L:NonIneq} has to be strengthened into the
following:

\setcounter{claim}{0}
\begin{claim}
The following relations hold:
 \[
 \gh_1\nleq\gx_1\jj\gz_1;\quad\gh_2\nleq\gx_0\jj\gz_3;\quad
 \gh_3\nleq\gx_3\jj\gz_2;\quad\gh_0\nleq\gx_2\jj\gz_0.
 \]
\end{claim}

\begin{cproof}
These relations follow, respectively, from the relations
 \[
 4\in\gh_1\setminus(\gx_1\jj\gz_1);\quad
 5\in\gh_2\setminus(\gx_0\jj\gz_3);\quad
 6\in\gh_3\setminus(\gx_3\jj\gz_2);\quad
 7\in\gh_0\setminus(\gx_2\jj\gz_0).
 \]
\end{cproof}

The proof of Theorem~\ref{T:NonAmalg2} proceeds then as follows.
We choose the elements
$0_K$, $1_K\in K$, $0_{K_i}=f_i(0_K)$, $1_{K_i}=f_i(1_K)$,
$0_{L_j}=g_{ij}(0_{K_i})$, for $i\ne j$, $i$, $j<3$, and $0_P$,
$1_P\in P$ in the same way.

Now, because each $S_i$ has almost permutable congruences, there
are elements $x_i\in S_i$ such that $0_{K_i}<x_i<1_{K_i}$ and
 \[
 \gF_{K_i}(0_{K_i},x_i)\in\set{\ga_i,\gb_i},\quad
 \gF_{K_i}(x_i,1_{K_i})\in\set{\ga_i,\gb_i}
 \]
and $\gF_{K_i}(0_{K_i},x_i)\ne\gF_{K_i}(x_i,1_{K_i})$. 

The rest consists of considering various combinations for the elements
$\gF_{K_i}(0_{K_i},x_i)$ and $\gF_{K_i}(x_i,1_{K_i})$. The
proof of Theorem~\ref{T:NonAmalg} works in the case
$\gF_{K_i}(0_{K_i},x_i)=\ga_i$ for $i<3$, and also in the case
$\gF_{K_i}(0_{K_i},x_i)=\gb_i$, $\gF_{K_i}(x_i,1_{K_i})=\ga_i$.
In the latter case we only need to dualize the original proof, which
leads to the same contradiction.

The next case, $\gF_{K_0}(0_{K_0},x_0)=\gb_0$,
$\gF_{K_1}(0_{K_1},x_1)=\ga_1$, and
$\gF_{K_2}(0_{K_2},x_2)=\gb_2$ consists of exchanging $\ga_0$ and
$\gb_0$, $\ga_2$ and $\gb_2$ in the original proof. This
case leads to the inequality $\gh_2\leq\gx_0\jj\gz_3$ which
contradicts the second case of Claim~1. The case
$\gF_{K_0}(0_{K_0},x_0)=\ga_0$,
$\gF_{K_1}(0_{K_1},x_1)=\gb_1$,
$\gF_{K_2}(0_{K_2},x_2)=\ga_2$ is dual to the former one and
leads to the same contradiction.

The case $\gF_{K_0}(0_{K_0},x_0)=\gb_0$,
$\gF_{K_1}(0_{K_1},x_1)=\gb_1$, and
$\gF_{K_2}(0_{K_2},x_2)=\ga_2$ leads to the inequality
$\gh_3\leq\gx_3\jj\gz_2$, contradicting the third case of Claim~1.
The case $\gF_{K_0}(0_{K_0},x_0)=\ga_0$,
$\gF_{K_1}(0_{K_1},x_1)=\ga_1$, and
$\gF_{K_2}(0_{K_2},x_2)=\gb_2$ is dual. 

Finally, the case $\gF_{K_0}(0_{K_0},x_0)=\ga_0$,
$\gF_{K_1}(0_{K_1},x_1)=\gb_1$, and
$\gF_{K_2}(0_{K_2},x_2)=\gb_2$ leads to the inequality
$\gh_0\leq\gx_2\jj\gz_0$, contradicting the last case of Claim~1,
while the remaining case $\gF_{K_0}(0_{K_0},x_0)=\gb_0$,
$\gF_{K_1}(0_{K_1},x_1)=\ga_1$,
$\gF_{K_2}(0_{K_2},x_2)=\ga_2$ is dual to the previous one.
\end{proof}

\section{No functorial solution of the Congruence Lattice
Problem}\label{S:NoFunc}

We shall give in this section a very elementary
diagram of finite Boolean semilattices and
\jz-homomorphisms, that cannot be lifted,
in an isomorphism-preserving fashion, by lattices and
lattice homomorphisms. This diagram is displayed on
Figure~\ref{Fi:DiagSemil2}:

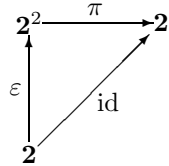
\begin{figure}[htb]
\begin{picture}(50,50)(0,5)

\put(0,0){\makebox(0,0){$\two$}}
\put(0,50){\makebox(0,0){$\two^2$}}
\put(50,50){\makebox(0,0){$\two$}}

\put(0,5){\vector(0,1){40}}
\put(5,50){\vector(1,0){42}}
\put(3.5,3.5){\vector(1,1){42}}

\put(-5,25){\makebox(0,0){$\ge$}}
\put(25,55){\makebox(0,0){$\gp$}}
\put(30,20){\makebox(0,0){$\id$}}

\end{picture}
\caption{Triangular semilattice diagram}
\label{Fi:DiagSemil2}
\end{figure}

The semilattice maps $\ge$ and $\gp$ are defined as
follows: $\ge(x)=\vv<x,x>$ and
$\gp(\vv<x,y>)=x\jj y$, for all $x$, $y<2$. Thus it is obvious
that the diagram of Figure~\ref{Fi:DiagSemil2} is commutative.

The proof of the following fact is so simple-minded that it
hardly deserves to be called a theorem. However, it implies
immediately that the Congruence Lattice Problem does not have a
functorial solution from \jz-semilattices and \jz-homomorphisms,
to lattices and lattice homomorphisms, see
Corollary~\ref{C:nofunct}.

\begin{theorem}\label{T:nofunct}
There is no lifting, with respect to the $\Conc$ functor, in
the category of lattices, of the semilattice diagram displayed on
Figure~~\tup{\ref{Fi:DiagSemil2}}, that sends the identity to
an isomorphism.
\end{theorem}

\begin{proof}
Assume, to the contrary, that the diagram can be lifted,
by a lattice diagram of the format displayed on
Figure~\ref{Fi:DiagLatt2}, with $f$ surjective.

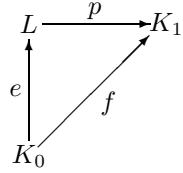
\begin{figure}[htb]
\begin{picture}(50,50)(0,5)

\put(0,0){\makebox(0,0){$K_0$}}
\put(0,50){\makebox(0,0){$L$}}
\put(52,50){\makebox(0,0){$K_1$}}

\put(0,6){\vector(0,1){38}}
\put(5,50){\vector(1,0){41}}
\put(3.5,3.5){\vector(1,1){42}}

\put(-5,25){\makebox(0,0){$e$}}
\put(25,55){\makebox(0,0){$p$}}
\put(30,20){\makebox(0,0){$f$}}

\end{picture}
\caption{Triangular lattice diagram}
\label{Fi:DiagLatt2}
\end{figure}

In particular, $p\circ e=f$ is surjective, thus $p$
is surjective. On the other
hand, $\Conc p$ is isomorphic (as a semilattice
homomorphism) to $\gp$, and $\gp$ separates $0$, thus $p$
is one-to-one. Therefore, $p$ is an isomorphism, which is
impossible since $\gp\iso\Conc p$ and $\gp$ is not an
isomorphism.
\end{proof}

\begin{corollary}\label{C:nofunct}
There is no quasi-functor $\FF$, from \jz-semilattices and
\jz-homomorphisms to lattices and lattice homomorphisms, such
that $\Conc\FF(S)\iso S$ and $\Conc\FF(f)\iso f$, for all finite
Boolean semilattices $S$ and $T$ and all \jz-homomorphisms
$f\colon S\to T$.
\end{corollary}

(The definition of a \emph{quasi-functor} is similar to the
definition of a functor, except that the image of an identity is
not required to be an identity.)

\begin{proof}
If $\gn\colon\two\to\two$ is the identity, then
$\gn\circ\gn=\gn$. Hence, if $f=\FF(\gn)$, then $f\circ f=f$.
However, by assumption on $\FF$, the relation $\Conc f\iso\gn$
holds. Since $\gn$ separates zero, $f$ is one-to-one, thus,
since $f$ is idempotent, $f=\id_{\FF(\two)}$. By
Theorem~\ref{T:nofunct}, this is impossible.
\end{proof}

Of course, the map $\gp$ is not one-to-one, and, in
particular, the proof of
Theorem~\ref{T:nofunct} does not imply the non-existence of
a functor from \jz-semilattices \emph{with \jz-embeddings}
to lattices and lattice homomorphisms, that lifts the
$\Conc$ functor. In fact, the diagram of
Figure~\ref{Fi:DiagSemil2} has a lifting by finite lattices
and lattice homomorphisms, as follows. In
Figure~\ref{Fi:DiagLatt2}, define $K_0=\two$, $K_1=M_3$ (the
five-element modular non-distributive lattice), $L=\two^2$;
let $e$ and $f$ be the $\set{0,1}$-preserving maps, and let
$p$ be any embedding from $\two^2$ into $M_3$.

\section{Open problems}

By a result of P. Pudl\'ak, see Fact 4, page 100 in \cite{Pudl},
every distributive \jz-semilattice is the direct union of all its
finite distributive \jz-subsemilattices. Therefore, in view of
the negative results of Section~\ref{S:NoFunc}, a positive
solution to the following Problem~\ref{Pb:Funct} would be about
the best possible solution to the Congruence Lattice Problem:

\begin{problem}\label{Pb:Funct}
Does there exist a functor $\FF$, from finite distributive
\jz-semilattices and their embeddings to lattices and their
embeddings, such that the functor $\Conc\circ\FF$ is naturally
equivalent to the identity?
\end{problem}

A related open problem is the following:

\begin{problem}\label{Pb:exlift}
Does every finite diagram (indexed by a poset) of finite
distributive \jz-semilattices have a lifting, with respect to
the $\Conc$ functor, by a diagram of lattices?
\end{problem}

We have seen in Section~\ref{S:liftDc} that the diagram $\Dc$ can
be lifted with respect to the $\Conc$ functor. However, we do not
even know the general answer to the following problem, thus
illustrating the level of our ignorance about
Problem~\ref{Pb:exlift}:

\begin{problem}\label{Pb:LiftCube}
Let $\E D$ be a cube of finite distributive \jz-semilattices. Is
it decidable whether $\E D$ admits a lifting, with respect to the
$\Conc$ functor, by a cube of lattices (resp., lattices with
permutable congruences)?
\end{problem}

We do not even know the answer to Problem~\ref{Pb:LiftCube} in the
particular case where $\E D=\Dac$.

\begin{problem}\label{Pb:PermCon}
Which algebraic distributive lattices are isomorphic to $\Con L$,
for some lattice $L$ with permutable congruences?
\end{problem}

A first approach to Problem~\ref{Pb:PermCon} might be provided by
E. T. Schmidt's well-known sufficient condition, for a given
algebraic distributive lattice, to be isomorphic to the congruence
lattice of a lattice, see \cite{Schm68}.
By using the amalgamation technique of \cite{GLWe} in a
ring-theoretical context, the second author proved that every
distributive \jz-semilattice of cardinality at most $\aleph_1$
is isomorphic to $\Conc L$ for some sectionally complemented
modular $L$, see \cite{Wehr3}. Since every sectionally
complemented lattice has permutable congruences, this provides a
strong positive answer to Problem~\ref{Pb:PermCon} for algebraic
distributive lattices with at most $\aleph_1$ compact elements.

\begin{problem}\label{Pb:NDvarComb}
Let $\VV$ be a non-distributive variety of lattices. Does there
exist a $\two^3$-diagram $\E D$ of lattices and lattice
embeddings in $\VV$ such that the image of $\E D$ under $\Conc$
cannot be lifted by lattices with almost permutable congruences?
\end{problem}

As follows from Corollary~\ref{C:NonReprAC}, if $\VV$ is a
non-distributive variety of lattices and if $F$ is a free
lattice in $\VV$ on at least $\aleph_2$ generators, then there
exists no lattice $L$ with almost permutable congruences such
that $\Conc L\iso\Conc F$. So, Problem~\ref{Pb:NDvarComb} asks
for a combinatorial analogue of Corollary~\ref{C:NonReprAC}.

\end{document}